# INFINITE FAMILIES OF VERY EXOTIC SPHERES WITH FREE $S^1$- AND $S^3$-ACTIONS

TILMAN BAUER AND J.D. QUIGLEY

Abstract. There are two kinds of exotic spheres: bp spheres, which bound parallelizable manifolds, and non-bp spheres, or very exotic spheres, which do not. In the 1960s, W.-C. Hsiang showed that in each dimension where bp spheres exist, there is at least one which admits infinitely many inequivalent smooth free $S^1$-actions, and in each dimension congruent to 3 modulo 4, there is at least one bp sphere which admits infinitely many inequivalent smooth free $S^3$-actions. On the other hand, for each fixed prime $p$, smooth free $S^1$- and $S^3$-actions have only been recorded to exist for finitely many very exotic spheres with nontrivial $p$-local Kervaire–Milnor invariant, all in dimension less than approximately $p^3$. In this paper, we use topological modular forms to detect smooth free $S^1$- and $S^3$-actions on infinite families of very exotic spheres with nontrivial 2- and 3-local Kervaire–Milnor invariants.

## 1. Introduction

An *exotic $n$-sphere* is a smooth $n$-dimensional manifold which is homeomorphic, but not diffeomorphic, to the $n$-sphere $S^n \subseteq \mathbb{R}^{n+1}$ equipped with its standard smooth structure. Wang and Xu [WX17] conjectured that exotic $n$-spheres exist for all $n \geq 5$ except for $n \in \{5, 6, 12, 56, 61\}$. This conjecture is true for $n$ odd by [WX17, HHR16, Bro69, KM63], for $n \leq 139$ by [BHHM20], and for at least half of all even $n$ by [BMQ23]. Although many exotic spheres are known to exist, relatively little is known about their geometry, and in particular their symmetry.

**Question 1.1** (Hsiang–Hsiang [HH68])**.** Does every exotic sphere admit a smooth faithful $S^1$-action?

All exotic spheres which bound parallellizable manifolds ("bp spheres") admit smooth faithful $S^1$-actions via equivariant plumbing or analysis of Brieskorn representations (cf. [HH67]), so we restrict attention to exotic spheres which do not bound parallelizable manifolds. We will refer to these as "non-bp spheres", or, following Schultz [Sch85], as "very exotic spheres". Bredon [Bre67] and Schultz [Sch75] showed that certain non-bp spheres of dimensions $(8k+1)$- and $(8k+2)$ (those with nontrivial $\alpha$-invariant) admit smooth faithful $S^1$-actions. The second author's results with Botvinnik in [BQ23] imply that infinitely many other non-bp spheres admit smooth faithful $S^1$-actions, although these only account for a small proportion of all non-bp spheres.

In this paper, we will focus on the following variant of the Hsiangs' question:

**Question 1.2.** Which exotic spheres admit a smooth *free* $S^1$-action?

Since $S^1$ cannot act continously and freely on an even-dimensional sphere, the question is only relevant for odd-dimensional exotic spheres. In [Hsi66], W.-C. Hsiang showed that bp spheres in every odd dimension where they exist admit smooth free $S^1$-actions, so again, we restrict attention to non-bp spheres. R. Lee [Lee68] showed that the exotic $(8k+1)$-dimensional non-bp spheres with nontrivial $\alpha$-invariant do *not* admit free $S^1$-actions. Moreover, none of the actions obtained using the results of [BQ23] will be free.

To state some positive results, recall that up to summing with bp spheres, non-bp spheres are classified by their Kervaire–Milnor invariant, i.e., the image of their Pontryagin–Thom invariant in the cokernel of the stable J-homomorphism [KM63]. To the authors' knowledge, the only non-bp spheres which are known to admit free $S^1$-actions are as follows:

(1) For $p = 2$, the non-bp spheres in dimensions 9, 15, 21, and 23 with 2-local Pontryagin–Thom invariants $\nu^3$, $\kappa\eta$, $\sigma^3$, and $\bar{\kappa}\nu$, respectively;
(2) For $p$ an odd prime, the non-bp spheres in dimensions less than $2p^3 + 2p^2 - 4p - 3$ with $p$-local Pontryagin–Thom invariants $\alpha_1\beta_1^r$ $(r < p)$, $\alpha_1\beta_1^r\beta_s$ $(r + s \neq p)$, $\alpha_1\epsilon_i$ $(1 \leq i < p - 2)$, $\beta_1\epsilon'$, and $p\varphi$.

These were all proven to exist in unpublished work of Brumfiel [Bru69, Props. 8.14, 8.15], to which we refer for naming conventions. The action on the non-bp 9-sphere with 2-local Pontryagin–Thom invariant $\nu^3$ and the non-bp 13-sphere with 3-local Pontryagin–Thom invariant $\alpha_1\beta_1$ was published in [Bru71, Thm. 1.10], albeit citing the aforementioned unpublished computations. A fully published proof of the action on the 9-sphere was given by D.M. James in [Jam80] using different methods. The action on the non-bp 15-sphere was proven by Bredon [Bre67], also using different methods from Brumfiel, and this action was extended to a smooth free $S^3$-action using Brumfiel's techniques in recent work of Basu and Kasilingam [BK22]. In summary, only finitely many non-bp spheres with nontrivial $p$-local Pontryagin–Thom invariant are known to admit smooth free $S^1$-actions for each prime $p$.

In this paper, we identify infinite families of non-bp spheres with nontrivial 2- and 3-local Pontryagin–Thom invariants which admit smooth free $S^1$- and $S^3$-actions.

**Theorem A.**

(1) For each $k \in \mathbb{N}$, there exists a $(15 + 96k)$-dimensional non-bp sphere with nontrivial 2-local Pontryagin–Thom invariant and a $(13 + 72k)$-dimensional non-bp sphere with nontrivial 3-local Pontryagin–Thom invariant admitting a smooth free $S^1$-action.
(2) For each $k \in \mathbb{N}$, there exists a $(35 + 96k)$-dimensional, a $(15 + 96k)$- or $(11 + 96k)$-dimensional, and a $(59+96k)$- or $(55+96k)$-dimensional non-bp sphere with nontrivial 2-local Pontryagin–Thom invariant admitting a smooth free $S^3$-action.

To prove the theorem, the first step is translating from geometric topology to stable homotopy theory, which is accomplished using the classical ideas of Brumfiel mentioned above [Bru69, Bru71], as well as a construction of Bredon [Bre67].

Roughly speaking, Brumfiel's work allows one to detect smooth free $S^1$-actions on non-bp $(2n+1)$-spheres by producing nontrivial elements in $\pi^0_s(\mathbb{C}P^n)$, the zeroth stable cohomotopy group of the $n$th complex projective space, which do not extend to $\pi^0_s(\mathbb{C}P^{n+k})$ for some $k > 0$. The same ideas can be adapted to study free $S^3$-actions on non-bp $(4n + 3)$-spheres using the stable cohomotopy of quaternionic projective spaces. We recall Brumfiel's work on $S^1$-actions and work out the analogous results for $S^3$-actions in Section 2, with the reduction from free actions to stable cohomotopy appearing in Section 2.1, Section 2.2, and Section 2.3. A summary is given in Section 2.5.

Bredon's construction and applications to $S^1$- and $S^3$-actions are recalled in Section 2.6. Roughly speaking, Bredon defined a pairing which, given a homotopy sphere and an element in the stable stems, produces a new homotopy sphere. Applying this pairing with elements from the stable 1- and 3-stems provides a method for detecting free $S^1$- and $S^3$-actions, respectively, in terms of relations in the stable stems.

The second step is solving the relevant problems in stable homotopy theory. For Bredon's construction, we simply check the stable stems for instances of the relevant relations. All of these relations are detected using *connective topological modular forms*, a generalized cohomology theory which detects infinite families of nontrivial elements in the stable stems. We refer to [DFHH14] for an introduction to *tmf*.

For Brumfiel's approach, the problem is reduced to finding nontrivial torsion elements in $\pi_s^0(\mathbb{C}P^n)$ and $\pi_s^0(\mathbb{H}P^n)$ for various $n$ after localizing at 2 or 3. To do so, we use the *tmf*-Hurewicz homomorphisms

$$\pi_s^0(\mathbb{C}P^n) \to tmf^0(\mathbb{C}P^n) \quad \text{and} \quad \pi_s^0(\mathbb{H}P^n) \to tmf^0(\mathbb{H}P^n).$$

For our purposes, the essential points regarding *tmf* are that the Hurewicz homomorphisms for spheres $\pi_s^0(S^n) \to tmf^0(S^n)$ are completely understood [BMQ23, BS23], the *tmf*-cohomology of $\mathbb{C}P^n$ is completely computable after localizing at 3, and the *tmf*-cohomology of $\mathbb{H}P^n$ is computable after localizing at 2 or 3. The analysis of these *tmf*-cohomology groups is the technical core of the paper, and is carried out in Section 4.1, Section 4.2, and Section 5.2.

These computations and various cofiber sequences involving complex and quaternionic projective spaces allow us to lift nontrivial elements in $tmf^0\mathbb{C}P^n$ and $tmf^0\mathbb{H}P^n$ to $\pi_s^0\mathbb{C}P^n$ and $\pi_s^0\mathbb{H}P^n$. The following theorem combines the results of Section 4.3 and Section 5.3.

**Theorem B.** For each $j \geq 0$:

(1) There exist nontrivial elements in $(\pi_s^0\mathbb{H}P^{2+24j})_{(2)}$ and $(\pi_s^0\mathbb{H}P^{13+24j})_{(2)}$ with nontrivial Hurewicz image in $(tmf^0\mathbb{H}P^{2+24j})_{(2)}$ and $(tmf^0\mathbb{H}P^{13+24j})_{(2)}$, respectively.

(2) There exists a nontrivial element in $(\pi_s^0\mathbb{C}P^{6+36j})_{(3)}$ with nontrivial Hurewicz image in $(tmf^0\mathbb{C}P^{6+36j})_{(3)}$.

The third and final step is translating back from stable homotopy theory to geometric topology. In Bredon's work, this is automatic, but in Brumfiel's work, more care is required. Using the surgery exact sequence, each element in part (1) of the previous theorem detects a nontrivial homotopy smoothing of $\mathbb{H}P^{2+24j}$ and $\mathbb{H}P^{13+24j}$, and each element in part (2) detects a nontrivial homotopy smoothing of $\mathbb{C}P^{6+36j}$. These homotopy smoothings correspond to smooth free $S^3$- and $S^1$-actions, respectively, but as we explain in Section 2.4, determining the precise sphere on which we obtain an action requires some additional analysis involving relations in the stable homotopy groups of spheres. This is carried out in Section 4.4 and Section 5.4. The ambiguity in part (1) of Theorem A stems from our incomplete understanding of certain explicit products in the stable homotopy groups of spheres.

1.1. **Acknowledgements.** The authors thanks Gabriel Angelini-Knoll, Boris Botvinnik, Reinhard Schultz, and Oliver Wang for helpful discussions. The second author was supported by NSF grants DMS-2414922 and DMS-2441241.

## 2. Reductions to stable homotopy theory

In this section, we recall Brumfiel and Bredon's work on $S^1$- and $S^3$-actions on homotopy spheres from [Bre67, Bru69, Bru71], the former of which builds on ideas of W.C. Hsiang from [Hsi66]. We begin with Brumfiel's work, and since the primary reference [Bru69] is unpublished, we recall some details. In Section 2.1, we explain the passage from free $S^1$- and $S^3$-actions on homotopy spheres to homotopy smoothings of complex and quaternionic projective spaces. In Section 2.2, we describe the passage from homotopy smoothings to $F/O$-bundles. In Section 2.3, we describe the reduction from $F/O$-bundles to stable cohomotopy of complex and quaternionic projective spaces. In Section 2.4, we discuss how to determine the homotopy sphere on which one has detected an action through the above procedure. We summarize the results in Section 2.5. We then give a rapid recollection of key results from Bredon's work [Bre67], which allows one to detect free $S^1$- and $S^3$-actions on homotopy spheres using relations in the stable homotopy groups of spheres in Section 2.6.

2.1. **From free actions to homotopy smoothings.** Suppose $\Sigma^{2n+1}$ is a homotopy $(2n+1)$-sphere with homotopy smoothing

$$\tilde{h} : \Sigma^{2n+1} \xrightarrow{\simeq} S^{2n+1}$$

which is equipped with a smooth free $S^1$-action

$$T : S^1 \times \Sigma^{2n+1} \to \Sigma^{2n+1}.$$

Then we have a commutative diagram

$$\begin{array}{ccc} \Sigma^{2n+1} & \xrightarrow{\tilde{h}} & S^{2n+1} \\ \downarrow & & \downarrow \\ P^n & \xrightarrow{h} & \mathbb{C}P^n, \end{array}$$

where $P^n := \Sigma^{2n+1}/S^1$ and $h : P^n \to \mathbb{C}P^n$ classifies the principal $S^1$-bundle over $P^n$ determined by the $S^1$-action $T$. The map $h$ is a homotopy equivalence, so $P^n$ is a homotopy smoothing of $\mathbb{C}P^n$. This gives a function from the set of free circle actions on homotopy $(2n+1)$-spheres to the set of homotopy smoothings of $\mathbb{C}P^n$. In fact, we have the following:

**Proposition 2.1** ([Bru69, p. 3] or [Hsi66, Sec. 2])**.** For $n \geq 3$, the function defined above provides a bijection between the set of equivariant diffeomorphism classes of smooth free $S^1$-actions on homotopy $(2n+1)$-spheres and the set of equivalence classes of homotopy smoothings of $\mathbb{C}P^n$.

**Notation 2.2.** If $M$ is a smooth manifold, we write $hS(M)$ for the surgery structure set of $M$, i.e., the set of isotopy classes of homotopy smoothings of $M$. We regard $hS(M)$ as a pointed set with basepoint $M \xrightarrow{id} M$.

Completely analogous arguments prove the following analogue for $S^3$-actions.

**Proposition 2.3.** For $n \geq 2$, there is a bijection between the set of equivariant diffeomorphism classes of smooth free $S^3$-actions on homotopy $(4n+3)$-spheres and $hS(\mathbb{H}P^n)$.

2.2. **From homotopy smoothings to $F/O$-bundles.** Let $M^k$ be a smooth $k$-manifold and let $F/O := \mathrm{fib}(BO \to BF)$. An *$F/O$-bundle* is a map $g : M^k \to F/O$. Geometrically, an $F/O$-bundle is a stable vector bundle $\xi(g) \to M^k$ together with a fiber homotopy trivialization $G : S(\xi) \to S^{\dim(\xi)}$ of its sphere bundle $S(\xi) \to M^k$. Making $G$ transverse regular to a point $p \in S^N$ gives a framed submanifold

$$L^k \times \mathbb{R}^N \subseteq S(\xi)$$

such that

- $G|_{L^k \times \mathbb{R}^N} : L^k \times \mathbb{R}^N \to S^N - p \cong \mathbb{R}^N$ is the projection onto $\mathbb{R}^N$,
- $G(S(\xi) - L^k \times \mathbb{R}^N) = p$, and
- $\pi : L^k \to M^k$ is of degree one.

Given a framed submanifold $L^k \times \mathbb{R}^N \subseteq S(\xi)$, we can do framed surgery to try to make $\pi : L^k \to M^k$ into a homotopy equivalence without changing the associated $F/O$-bundle over $M^k$. The only obstruction is the middle-dimensional surgery obstruction, which for simply connected $M$ lies in the group $P_k = L_k(\mathbb{Z})$, which is $\mathbb{Z}$ if $k \equiv 0 \mod 4$, $\mathbb{Z}/2$ if $k \equiv 2 \mod 4$, and zero otherwise.

Conversely, given a homotopy equivalence $f : L^k \to M^k$ with homotopy inverse $\bar{f}$, we can cover $f$ up to homotopy by an embedding of $L^k$ into the total space of $\xi = \bar{f}^*(\tau_L) - \tau_M$ over $M^k$. Since $L^k$ has trivial normal bundle in $\xi$, we can choose a framing $L^k \times \mathbb{R}^N \subseteq \xi$ to define a map

$$\theta(L^k, f) : M^k \to F/O.$$

One can show that the homotopy class of $\theta(L^k, f)$ is independent of choices.

Below, we will be interested in the homotopy class of $\theta(L^k, f)$. Since $F/O$ is a path-connected H-space, unpointed homotopy classes of maps from $M^k$ into $F/O$ are in bijection with pointed homotopy classes of maps from $M^k$ into $F/O$. It will be slightly more convenient below to work with pointed homotopy classes, so we state the following proposition using pointed maps.

**Proposition 2.4.** Let $M^k$ be a simply connected smooth $k$-manifold. The surgery exact sequence (cf. [Ran02, 1.18]) is an exact sequence of pointed sets

$$P_{k+1} \to hS(M^k) \xrightarrow{\theta} [M^k, F/O] \xrightarrow{s} P_k,$$

where $[M^k, F/O]$ denotes the pointed set of pointed homotopy classes of maps from $M^k$ to $F/O$. In particular, if $k$ is even, then $\theta$ is injective.

Specializing to the cases $M^k = \mathbb{C}P^n$ and $M^k = \mathbb{H}P^n$, we obtain the following:

**Corollary 2.5.**

(1) For all $n \geq 3$, there is an exact sequence of pointed sets

$$0 \to hS(\mathbb{C}P^n) \xrightarrow{\theta} [\mathbb{C}P^n, F/O] \xrightarrow{s} P_{2n}.$$

If $n$ is even, then $P_{2n} = \mathbb{Z}$, and if $n$ is odd, then $P_{2n} = \mathbb{Z}/2$.

(2) For all $n \geq 2$, there is an exact sequence of pointed sets

$$0 \to hS(\mathbb{H}P^n) \xrightarrow{\theta} [\mathbb{H}P^n, F/O] \xrightarrow{s} \mathbb{Z}.$$

2.3. **From $F/O$-bundles to stable cohomotopy.** To translate from $F/O$-bundles to stable cohomotopy, we will use the sequence of fibrations

$$SO \to SF \to F/O \to BSO \to BSF. \tag{1}$$

The following result combines [Bru69, Cors. 2.1 and 2.3]. Since we will need to modify the result in the quaternionic case, we sketch the proof from Brumfiel's work.

**Proposition 2.6** ([Bru69]). There is a split short exact sequence

$$0 \to [\mathbb{C}P^n, F] \to [\mathbb{C}P^n, F/O] \to \bigoplus_{i=1}^{\lfloor n/2 \rfloor} \mathbb{Z} \to 0.$$

In particular, there are isomorphisms

$$[\mathbb{C}P^n, F/O] \cong \bigoplus_{i=1}^{\lfloor n/2 \rfloor} \mathbb{Z} \oplus [\mathbb{C}P^n, F] \quad \text{and} \quad [\mathbb{C}P^n, F/O]_{\mathrm{tor}} \cong [\mathbb{C}P^n, F] \cong \pi_s^0(\mathbb{C}P^n).$$

*Proof.* Applying $[\mathbb{C}P^n, -]$ to the sequence (1) yields an exact sequence

$$[\mathbb{C}P^n, SO] \to [\mathbb{C}P^n, F] \to [\mathbb{C}P^n, F/O] \to [\mathbb{C}P^n, BSO] \to [\mathbb{C}P^n, BF].$$

An elementary Atiyah–Hirzebruch spectral sequence computation shows that

$$[\mathbb{C}P^n, SO] \cong KO^{-1}\mathbb{C}P^n = 0,$$

yielding injectivity of the map $[\mathbb{C}P^n, F] \to [\mathbb{C}P^n, F/O]$. Moreover, Adams and Walker [AW65] showed that

$$[\mathbb{C}P^n, BSO] \cong KO^0\mathbb{C}P^n \cong \begin{cases} \bigoplus_{i=1}^{\lfloor n/2 \rfloor} \mathbb{Z} & \text{if } n \equiv 0, 2, 3 \mod 4, \\ \left(\bigoplus_{i=1}^{\lfloor n/2 \rfloor} \mathbb{Z}\right) \oplus \mathbb{Z}/2 & \text{if } n \equiv 1 \mod 4. \end{cases}$$

Adams and Walker also proved that the kernel of $[\mathbb{C}P^n, BSO] \to [\mathbb{C}P^n, BF]$ is a free abelian group of maximal rank. The proposition follows. □

We have the following quaternionic analogue.

**Proposition 2.7.** There is an exact sequence

$$\bigoplus_{i=1}^{\lfloor n/2 \rfloor} \mathbb{Z}/2 \to [\mathbb{H}P^n, F] \to [\mathbb{H}P^n, F/O] \to \bigoplus_{i=1}^{n} \mathbb{Z} \to 0.$$

*Proof.* Applying $[\mathbb{H}P^n, -]$ to the sequence (1) yields an exact sequence

$$[\mathbb{H}P^n, SO] \to [\mathbb{H}P^n, F] \to [\mathbb{H}P^n, F/O] \to [\mathbb{H}P^n, BSO] \to [\mathbb{H}P^n, BF].$$

The Atiyah–Hirzebruch spectral sequence shows that

$$[\mathbb{H}P^n, SO] \cong KO^{-1}\mathbb{H}P^n = \bigoplus_{i=1}^{\lfloor n/2 \rfloor} \mathbb{Z}/2.$$

No hidden extensions are possible because every class in $KO^{-1}$ is divisible by $\eta$, which is simple 2-torsion. On the other hand, the quaternionic analogues of the results of Adams and Walker proven by Sigrist and Suter [SS74] imply that

$$[\mathbb{H}P^n, BSO] \cong KO^0\mathbb{H}P^n \cong \bigoplus_{i=1}^{n} \mathbb{Z}$$

for all $n$, and since $[\mathbb{H}P^n, BF]$ is finite, we see that the kernel of $[\mathbb{H}P^n, BSO] \to [\mathbb{H}P^n, BF]$ is free abelian of maximal rank. □

**Remark 2.8.** The image of $\bigoplus_{i=1}^{\lfloor n/2 \rfloor} \to [\mathbb{H}P^n, F]$ is the image of the J-homomorphism for $\mathbb{H}P^n$. Our examples below will never be in the image of this map, and so they always have nontrivial image in $[\mathbb{H}P^n, F/O]$.

**Definition 2.9.** The groups $\pi_s^0\mathbb{C}P^n := [\mathbb{C}P^n, F]$ and $\pi_s^0\mathbb{H}P^n := [\mathbb{H}P^n, F]$ are the *zeroth stable cohomotopy groups* of $\mathbb{C}P^n$ and $\mathbb{H}P^n$.

We will be interested in lifting elements from $[\mathbb{C}P^n, F] \subseteq [\mathbb{C}P^n, F/O]$ and $\mathrm{Im}([\mathbb{H}P^n, F] \to [\mathbb{H}P^n, F/O])$ to $hS(\mathbb{C}P^n)$ and $hS(\mathbb{H}P^n)$. To do so, we will use the following.

**Proposition 2.10.**

(1) For $n = 4k,\ 4k+1,\ 4k+2$, and conjecturally if $n = 4k+3 \neq 2^j - 1$, we have $s([\mathbb{C}P^n, F]) = 0$. Moreover, if $n = 4k+3$, then any odd torsion classes in $[\mathbb{C}P^n, F]$ has trivial image under $s$.

(2) For any $n \geq 1$, we have $s(\mathrm{Im}([\mathbb{H}P^n, F] \to [\mathbb{H}P^n, F/O])) = 0$.

*Proof.* For the $\mathbb{C}P^n$ statement, see [Bru71, Lems. I.5, I.6, Conj. I.7], noting that $\mathbb{C}P^n$ in our notation is equivalent to $\mathbb{C}P(n+1)_0 = \mathbb{C}P^{n+1} \setminus D^{2n+2}$ in Brumfiel's notation. The claim about odd torsion follows from the observation that the target $P_{2n}$ has no odd torsion.

For the $\mathbb{H}P^n$ statement, observe that the image of $[\mathbb{H}P^n, F]$ in $[\mathbb{H}P^n, F/O]$ is torsion, while the target of $s$ is torsion-free. □

**Remark 2.11.** Although we will not need it, we note for the sake of completeness that when $n = 4k+3$, Brumfiel provides a method (attributed to Sullivan) for computing $s$ on $[\mathbb{C}P^n, F]$ using the total Wu class; see [Bru69, Lem. 8.3].

2.4. **The map $\sigma$.** The results above show that, in favorable circumstances, we can use elements in $\pi_s^0\mathbb{C}P^n$ and $\pi_s^0\mathbb{H}P^n$ to detect smooth free $S^1$- and $S^3$-actions on homotopy spheres. Given such an element, Brumfiel's work also describes how to determine *which* homotopy $(2n+1)$- or $(4n+3)$-sphere the action is on. We will summarize the discussion for $\mathbb{C}P^n$; the discussion for $\mathbb{H}P^n$ is completely analogous.

Let $\mathbb{C}P_0^{n+1} := \mathbb{C}P^{n+1} \setminus D^{2n+2}$ so that we have a $D^2$-bundle $D^2 \to \mathbb{C}P_0^{n+1} \xrightarrow{H} \mathbb{C}P^n$ and maps

$$S^{2n+1} \cong \partial\mathbb{C}P_0^{n+1} \overset{i}{\hookrightarrow} \mathbb{C}P_0^{n+1} \xrightarrow{H} \mathbb{C}P^n.$$

These induce maps

$$[\mathbb{C}P^n, F/O] \xrightarrow{H^*,\simeq} [\mathbb{C}P_0^{n+1}, F/O] \xrightarrow{\theta,\simeq} hS(\mathbb{C}P_0^{n+1}) \xrightarrow{i^*} \Gamma_{2n+1},$$

where $\Gamma_{2n+1}$ is the group of $h$-cobordism classes of oriented homotopy $(2n+1)$-spheres [KM63] and

$$i^*(h : M^{2n+2} \xrightarrow{\simeq} \mathbb{C}P_0^{n+1}) := \partial M^{2n+2}.$$

We define

$$\sigma := i^*\theta H^* : [\mathbb{C}P^n, F/O] \to \Gamma_{2n+1}. \tag{2}$$

The following result of Brumfiel explains the significance of $\sigma$.

**Proposition 2.12** ([Bru69, Prop. 1.1])**.** Let $f : P^{2n} \xrightarrow{\simeq} \mathbb{C}P^n$ correspond to $T : S^1 \times \Sigma^{2n+1} \to \Sigma^{2n+1}$. Then $\Sigma^{2n+1} = \sigma\theta(P^{2n}, f)$.

The work of Kervaire and Milnor [KM63], Brumfiel [Bru68], Mahowald and Tangora [MT67], Barratt, Jones, and Mahowald [BJM84], and Lin, Wang, and Xu [LWX24] implies that $\Gamma_{2n+1}$, the target of $\sigma$, admits the following decomposition:

**Proposition 2.13.** If $2n+1 \neq 5, 13, 29, 61, 125$, then $\Gamma_{2n+1} \cong bP_{2n+2} \oplus \operatorname{coker} J_{2n+1}$, where $bP_{2n+2}$ is the subgroup of homotopy spheres which bound parallelizable manifolds and $J_{2n+1}$ is the stable J-homomorphism.

Brumfiel showed [Bru71, Lem. I.8] that $\sigma(\bigoplus \mathbb{Z}) \subseteq bP_{2n+2}$, so classes from the free summand can only detect free $S^1$-actions on bp spheres. As these are relatively well-understood [Hsi66], we will restrict our attention to computing $\sigma$ on classes in $[\mathbb{C}P^n, F/O]$ and $[\mathbb{H}P^n, F/O]$ which come from the torsion groups $\pi_s^0\mathbb{C}P^n$ and $\pi_s^0\mathbb{H}P^n$.

Let $\varphi : \Gamma_{2n+1} \to \operatorname{coker} J_{2n+1}$ denote the Kervaire–Milnor invariant, i.e., the composite of the Pontryagin–Thom invariant with the projection onto the cokernel of $J$. Except in the exceptional dimensions listed above (the dimensions one less than those in which a manifold of Kervaire invariant one exists), a homotopy sphere $\Sigma$ is non-bp if and only if $\varphi(\Sigma) \neq 0$.

Since we are interested in detecting free $S^1$- and $S^3$-actions on non-bp spheres, our strategy will be to find classes $x \in \pi_s^0 \mathbb{C}P^n$ and $x \in \pi_s^0 \mathbb{H}P^n$ whose images under the composites

$$\pi_s^0 \mathbb{C}P^n \cong [\mathbb{C}P^n, F] \hookrightarrow [\mathbb{C}P^n, F/O] \xrightarrow{\sigma} \Gamma_{2n+1} \xrightarrow{\varphi} \operatorname{coker} J_{2n+1}, \tag{3}$$

$$\pi_s^0 \mathbb{H}P^n \cong [\mathbb{H}P^n, F] \to [\mathbb{H}P^n, F/O] \xrightarrow{\sigma} \Gamma_{4n+3} \xrightarrow{\varphi} \operatorname{coker} J_{4n+3} \tag{4}$$

are nonzero.

**Lemma 2.14.**

(1) For $2n+1 \neq 5, 13, 29, 61, 125$, the following diagram commutes:

$$\begin{array}{ccccccc}
\pi_s^0 \mathbb{C}P^n & \xrightarrow{\cong} & [\mathbb{C}P^n, F] & \xrightarrow{\circ H} & & [S^{2n+1}, F] & \xrightarrow{\cong} & \pi^s_{2n+1} \\
 & & \downarrow q & & & \downarrow q & & \downarrow q \\
 & & [\mathbb{C}P^n, F/O] & \xrightarrow{\sigma} & \Gamma_{2n+1} \xrightarrow{\varphi} & [S^{2n+1}, F/O] & \xrightarrow{\cong} & \operatorname{coker} J_{2n+1}
\end{array}$$

Here, each map labeled '$q$' is induced by the quotient map $F \to F/O$, and $H$ is the attaching map of the $(2n+2)$-cell of $\mathbb{C}P^{n+1}$,

$$S^{2n+1} \xrightarrow{H} \mathbb{C}P^n \hookrightarrow \mathbb{C}P^{n+1}.$$

(2) For all $n$, the following diagram commutes:

$$\begin{array}{ccccccc}
\pi_s^0 \mathbb{H}P^n & \xrightarrow{\cong} & [\mathbb{H}P^n, F] & \xrightarrow{\circ H'} & & [S^{4n+3}, F] & \xrightarrow{\cong} & \pi^s_{4n+3} \\
 & & \downarrow q & & & \downarrow q & & \downarrow q \\
 & & [\mathbb{H}P^n, F/O] & \xrightarrow{\sigma} & \Gamma_{4n+3} \xrightarrow{\varphi} & [S^{4n+3}, F/O] & \xrightarrow{\cong} & \operatorname{coker} J_{4n+3},
\end{array}$$

where each $q$ is induced by the quotient map $F \to F/O$ and $H'$ is the attaching map of the $(4n+4)$-cell of $\mathbb{H}P^{n+1}$,

$$S^{4n+3} \xrightarrow{H'} \mathbb{H}P^n \hookrightarrow \mathbb{H}P^{n+1}.$$

*Proof.* This follows from the definition of $\sigma$. $\square$

**Notation 2.15.** The cellular filtrations of $\mathbb{C}P^n$ and $\mathbb{H}P^n$ give rise to (nonunitally) multiplicative Atiyah–Hirzebruch spectral sequences

$$E_1 = z(\pi_s^* S)[z]/(z^{n+1}) \Rightarrow \pi_s^* \mathbb{C}P^n, \quad E_1 = z(\pi_s^* S)[z]/(z^{n+1}) \Rightarrow \pi_s^* \mathbb{H}P^n.$$

Here $z$ denotes a generator in dimension $-2$ and $-4$, representing the bottom cell of $\mathbb{C}P^n$ and $\mathbb{H}P^n$, respectively.

The previous lemma and a straightforward diagram chase imply the following.

**Corollary 2.16.**

(1) Suppose $\alpha \neq 0 \in \pi_s^0 \mathbb{C}P^n$ is detected by a class $z^k x$ in the Atiyah–Hirzebruch spectral sequence converging to $\pi_s^* \mathbb{C}P^n$, and suppose further that $d_r(z^k x) = z^{n+1} y$ in the Atiyah–Hirzebruch spectral sequence converging to $\pi_s^* \mathbb{C}P^{n+1}$. Then $\varphi(\sigma(\alpha)) = q(y) \in \operatorname{coker} J_{2n+1}$. In particular, there exists a non-bp sphere with Pontryagin–Thom invariant $y$ which admits a free $S^1$-action.

(2) Suppose $\alpha \neq 0 \in \operatorname{coker}(KO^{-1}\mathbb{H}P^n \to \pi_s^0\mathbb{H}P^n)$ is detected by a class $z^kx$ in the Atiyah–Hirzebruch spectral sequence converging to $\pi_s^*\mathbb{H}P^n$, and suppose further that $d_r(z^k) = z^{n+1}y$ in the Atiyah–Hirzebruch spectral sequence converging to $\pi_s^*\mathbb{H}P^{n+1}$. Then $\varphi(\sigma(\alpha)) = q(y) \in \operatorname{coker} J_{4n+3}$. In particular, there exists a non-bp sphere with Pontryagin–Thom invariant $y$ which admits a free $S^3$-action.

2.5. **Summary of Brumfiel approach.** For future reference, we briefly summarize the discussion above.

2.5.1. *$S^1$-actions.* Equivariant diffeomorphism classes of free $S^1$-actions on homotopy $(2n+1)$-spheres are in bijection with $hS(\mathbb{C}P^n)$, the set of isotopy classes of homotopy smoothings of $\mathbb{C}P^n$. There is an exact sequence

$$0 \to hS(\mathbb{C}P^n) \to [\mathbb{C}P^n, F/O] \xrightarrow{s} P_{2n}$$

for $n \geq 3$, so we are interested in classes in $[\mathbb{C}P^n, F/O]$ which are in the kernel of $s$. There is also a map $\sigma : hS(\mathbb{C}P^n) \to \Gamma_{2n+1}$ which recovers the particular homotopy sphere on which we obtain a free $S^1$-action.

There is a splitting $[\mathbb{C}P^n, F/O] \cong [\mathbb{C}P^n, F] \oplus \bigoplus \mathbb{Z}$. We have $\sigma(\bigoplus \mathbb{Z}) \subseteq bP_{2n+2} \subseteq \Gamma_{2n+1}$. As free $S^1$-actions on bp-spheres are known to exist in all odd dimensions, we restrict attention to $[\mathbb{C}P^n, F] \cong \pi_s^0\mathbb{C}P^n$.

If $n \equiv 0, 1, 2 \mod 4$, then $s([\mathbb{C}P^n, F]) = 0$. If $n \equiv 3 \mod 4$, then $s = 0$ on odd-torsion classes.

Each nontrivial class $\alpha \in \pi_s^0\mathbb{C}P^n$ is detected by a class $z^kx$ for some $x \in \pi_{2k}^s$ and $1 \leq k \leq n$ which survives in the Atiyah–Hirzebruch spectral sequence for $\pi_s^*\mathbb{C}P^n$. If $\alpha$ lifts to a nontrivial class $\tilde{\alpha} \in hS(\mathbb{C}P^n)$, i.e., if $\alpha \in \ker(s)$, then we have $\sigma(\tilde{\alpha}) = \beta$, where $d_r(z^kx) = z^{n+1}\beta$ in the Atiyah–Hirzebruch spectral sequence converging to $\pi_s^*\mathbb{C}P^{n+1}$.

2.5.2. *$S^3$-actions.* Equivariant diffeomorphism classes of free $S^3$-actions on homotopy $(4n+3)$-spheres are in bijection with $hS(\mathbb{H}P^n)$, the set of isotopy classes of homotopy smoothings of $\mathbb{H}P^n$. There is an exact sequence

$$0 \to hS(\mathbb{H}P^n) \to [\mathbb{H}P^n, F/O] \xrightarrow{s} \mathbb{Z}$$

for $n \geq 2$, so we are interested in classes in $[\mathbb{H}P^n, F/O]$ which are in the kernel of $s$. There is also a map $\sigma : hS(\mathbb{H}P^n) \to \Gamma_{4n+3}$ which recovers the particular homotopy sphere on which $S^3$ acts.

There is an exact sequence

$$0 \to \bigoplus \mathbb{Z}/2 \to \pi_s^0\mathbb{H}P^n \to [\mathbb{H}P^n, F/O] \to \bigoplus \mathbb{Z} \to 0.$$

We have $\sigma(\bigoplus \mathbb{Z}) \subseteq bP_{4n+4}$. Since free $S^3$-actions on bp-spheres are known to exist in all dimensions congruent to 3 modulo 4, we restrict attention to nontrivial classes in the cokernel of $\bigoplus \mathbb{Z}/2 \to \pi_s^0\mathbb{H}P^n$.

For all of these classes, the surgery obstruction $s$ vanishes. Each nontrivial class $\alpha$ in the cokernel is detected by a class $z^kx$ for some $x \in \pi_{4k}^s$ and $1 \leq k \leq n$ which survives int he Atiyah–Hirzebruch spectral sequence for $\pi_s^*\mathbb{H}P^n$. Writing $\tilde{\alpha} \in hS(\mathbb{H}P^n)$ for a lift of $\alpha$ to $hS(\mathbb{C}P^n)$, we have $\sigma(\tilde{\alpha}) = \beta$, where $d_r(z^kx) = z^{n+1}\beta$ in the Atiyah–Hirzebruch spectral sequence converging to $\pi_s^*\mathbb{H}P^{n+1}$.

2.6. **Summary of Bredon approach.** In contrast with the Brumfiel approach summarized above, which requires analysis of the stable cohomotopy groups of various projective spaces, the work of Bredon [Bre67] allows us to detect free $S^1$- and $S^3$-actions on homotopy spheres directly from knowledge about the stable homotopy groups of spheres. In *loc. cit.*, Bredon defined a pairing

$$\rho_{n,k} : \Theta_n \otimes \pi_k^s \to \Theta_{n+k}$$

for $k < n-1$ as follows. Let $\Sigma$ be a homotopy $n$-sphere, $B \subset \Sigma$ an embedded $n$-disk, and let $M^k \subseteq S^{n+k}$ be a framed $k$-manifold corresponding to an element $\alpha \in \pi_k^s$ under the Pontryagin–Thom isomorphism. Then $\rho_{n,k}(\Sigma, \alpha)$ is the homotopy $(n+k)$-sphere obtained by removing a tubular neighborhood of $M$ from $S^{n+k}$ and replacing it with $M^k \times (\Sigma - B)$. We refer to [Bre67, Sec. 1] for details.

Bredon used this pairing to produce interesting group actions on homotopy spheres via the following construction [Bre67, Sec. 4]. Let $G$ be a compact Lie group acting smoothly on $S^{n+k}$ and equip $S^{n+k}$ with a $G$-invariant Riemannian metric. Let $x \in S^{n+k}$ and let $G(x)$ be a principal orbit of dimension $k$. Using the $G$-action, a choice of framing of $G(x)$ at $x$ extends to a canonical framing of the normal bundle of the entire orbit $G(x)$, and this framed $k$-manifold in $S^{n+k}$ represents an element of $\pi_{n+k}(S^k)$. If $k < n-1$, then this is an element in $\pi_{n+k}^s$. Let us denote this element by $\alpha$.

For $G = S^1$ and $G = S^3$, Bredon computed $\alpha$ for the standard free $G$-actions on $S^{n+1}$ and $S^{n+3}$ (see [Bre67, pp. 444-446]) (so $n$ is even if $G = S^1$ and $n \equiv 0 \mod 4$ if $G = S^3$).

(1) Let $G = S^1$, so $\alpha \in \pi_1^s \cong \mathbb{Z}/2$. If $n \equiv 1 \mod 4$, then $\alpha = 0 \in \pi_1^s$, while if $n \equiv 3 \mod 4$, then $\alpha = \eta \neq 0 \in \pi_1^s$.
(2) Let $G = S^3$, so $\alpha \in \pi_3^s \cong \mathbb{Z}/24$. Then $\alpha = \frac{n+1}{4} - 1 \in \mathbb{Z}/24$.

On the other hand, Bredon showed [Bre67, Cor. 2.2] that the diagram

$$\begin{array}{ccc} \Theta_n \otimes \pi_k^s & \xrightarrow{\rho_{n,k}} & \Theta_{n+k} \\ \downarrow{\scriptstyle p'\otimes 1} & & \downarrow{\scriptstyle p'} \\ \operatorname{coker}(J_n) \otimes \pi_k^s & \longrightarrow & \operatorname{coker}(J_{n+k}) \end{array}$$

commutes, where $p'$ is the Kervaire–Milnor invariant and the bottom arrow is induced by the composition product in $\pi_*^s$. Together with the observations above, we obtain the following homotopical criteria for detecting free $S^1$- and $S^3$-actions.

**Proposition 2.17.**

(1) Suppose $\Sigma$ is an exotic $(4\ell-1)$-sphere with $p'(\Sigma)$ divisible by $\eta$ in $\operatorname{coker}(J_{4\ell-1})_{(2)}$. Then $\Sigma$ admits a smooth free $S^1$-action.
(2) Suppose $\Sigma$ is an exotic $(4\ell-1)$-sphere with $p'(\Sigma)$ divisible by $[\ell]\nu$ in $\operatorname{coker}(J_{4\ell-1})_{(2)}$, where $[\ell]$ denotes the congruence class of $\ell$ modulo 8. Then $\Sigma$ admits a smooth free $S^3$-action.
(3) Suppose $\Sigma$ is an exotic $(4\ell-1)$-sphere with $\ell \equiv 0, 2 \mod 3$ and $p'(\Sigma)$ divisible by $\alpha_1$ in $\operatorname{coker}(J_{4\ell-1})_{(3)}$. Then $\Sigma$ admits a smooth free $S^3$-action.

**Remark 2.18.** Since we limit ourselves to using only topological modular forms for detection, we will only use the $p = 2$ parts of the previous proposition. An infinite family of examples for $p = 3$ can be deduced using the spectrum $j^2$ studied by Carrick and Davies [CD24], but we defer the use of $j^2$ and related cohomology theories to a future paper joint with the aforementioned authors.

## 3. Topological modular forms and the algebraic Atiyah–Hirzebruch spectral sequence

Brumfiel and Bredon's work reduce the detection of free $S^1$- and $S^3$-actions on exotic spheres to problems in stable homotopy theory. We will study both problems using the Hurewicz map for topological modular forms, $h : \mathbb{S} \to tmf$, which induces homomorphisms

$$h : \pi_n^s(X) \to tmf_n(X), \quad h : \pi_s^n(X) \to tmf^n(X)$$

from stable (co)homotopy groups to tmf (co)homology groups for all spectra $X$.

To apply Bredon's result (Proposition 2.17), we need to understand the actions of $\eta$, multiples of $\nu$, and $\alpha_1$ on the stable homotopy groups of spheres. Each of these elements is in the image of the *tmf*-Hurewicz homomorphism, so we may detect nontrivial actions on classes with nontrivial *tmf*-Hurewicz image by finding nontrivial actions in $tmf_*$. The *tmf*-Hurewicz image was computed after localizing at 2 by Behrens, Mahowald, and the second author in [BMQ23], and after localizing at 3 by Belmont and Shimomura in [BS23]. The image is trivial above degree zero after localizing at larger primes. These computations allow us to immediately detect infinite families of exotic spheres with free $S^1$- and $S^3$-actions.

Applying Brumfiel's result requires a bit more work. Our strategy for producing families of elements in the cohomotopy of $P^n$ (standing for either of $\mathbb{C}P^n$ or $\mathbb{H}P^n$) is to construct elements in $tmf^0(P^n)$ that are in the image of the Hurewicz homomorphism. Unlike cohomotopy, *tmf*-cohomology of projective spaces has two kinds of periodicity: one coming from the intrinsic 576-periodicity of *tmf* and the other from the Thom isomorphism derived from the string orientation of *tmf*, which implies that the *tmf*-cohomology of infinite projective space has a certain Thom class as a polynomial generator in it. Crucially, the image of the Hurewicz homomorphism in *tmf* is also 576-periodic. This means that a single class in $tmf^0(P^n)$ that is in the Hurewicz image gives rise to a periodic family of nonzero cohomotopy classes (as in Theorem B).

In this section, we briefly recall the Adams–Novikov spectral sequence for *tmf* and the algebraic Atiyah–Hirzebruch spectral sequence, which we will use extensively in the following computational sections.

The $MU$-based Adams–Novikov spectral sequence for $tmf \wedge X$ has signature

$$E_2 = \mathrm{Ext}^{**}_{(MU_*, MU_*MU)}(MU_*, (MU \wedge tmf)_*X) \Rightarrow tmf_*X.$$

By a change-of-rings isomorphism (cf. [DFHH14]), this $E_2$-term is isomorphic to the cohomology of the Weierstrass curve Hopf algebroid $(A, \Gamma)$, given by

$$A = \mathbb{Z}[a_1, a_2, a_3, a_4, a_6]; \quad \Gamma = A[r, s, t],$$

with coefficients in $A_*(X)$. Here, the homology theory $A_*(X)$ is constructed by $A_*(X) = (tmf \wedge X(4))_*(X)$ with coefficients $A_*(\mathrm{pt}) = A$, where $X(4)$ is a certain Thom spectrum first studied by Ravenel in [Rav84]; for details, we refer to [DFHH14]. From this point of view, the Adams–Novikov spectral sequence is sometimes called the descent spectral sequence.

Complete computations for $X$ a point can be found in [Bau08]. In this case, the $E_2$-term of the above spectral sequence is called the ring of *derived modular forms*. Its homological degree 0 part agrees with the ring of integral modular forms.

In order to understand the differentials in the *tmf*-based Atiyah–Hirzebruch spectral sequence for $P^n$, it is useful to consider the *algebraic Atiyah–Hirzebruch spectral sequence* instead. This spectral sequence has the signature

$$E_1 = \mathrm{Ext}^{**}_{(A,\Gamma)}(A, A)[z]/(z^{n+1}) \Longrightarrow \mathrm{Ext}^{**}_{(A,\Gamma)}(A, A_*(DP^n_+)).$$

It is a trigraded spectral sequence, where classes in $\mathrm{Ext}^{s,t}$ have tridegree $(s,t,0)$, $z$ has tridegree $(0,-4,-4)$ for $\mathbb{H}P^\infty$ and $(0,-2,-2)$ for $\mathbb{C}P^\infty$, and the degree of the differential $d^i$ is $(-1,0,i)$.

We consider this spectral sequence not only for $P^n$ but also, in an inductive argument, for the *truncated projective spaces*

$$P^n_k = P^n/P^{k-1} \cong (P^{n-k})^{kL},$$

where the last term denotes the Thom space of the $k$-fold sum of the tautological (complex or quaternion) line bundle on $P^{n-k}$.

The source and target of the algebraic Atiyah–Hirzebruch spectral sequence are $E_2$-terms of Adams–Novikov spectral sequences, and there is a commutative diagram of spectral sequences

$$\begin{array}{ccc} \mathrm{Ext}^{**}_{(A,\Gamma)}(A,A)[z] & \Longrightarrow & \mathrm{Ext}^{**}_{(A,\Gamma)}(A, A_*(D\mathbb{H}P^\infty_+)) \\ \Downarrow & & \Downarrow \\ \pi_* \mathit{tmf}[z] & \Longrightarrow & \mathit{tmf}_*(D\mathbb{H}P^\infty_+) \end{array}$$

where the bottom spectral sequence is the classical, *tmf*-based Atiyah–Hirzebruch spectral sequence.

Note that we choose to run the algebraic Atiyah-Hirzebruch spectral sequence at the same speed as the topological one, meaning that for $\mathbb{H}P^\infty$, $d^i = 0$ unless $4 \mid i$, and for $\mathbb{C}P^\infty$, $d^i = 0$ unless $2 \mid i$.

**Remark 3.1.** The spectrum $\Sigma^{-4}\mathit{tmf} \wedge \mathbb{H}P^n$ has been studied under different names. It is equivalent to the spectrum $\mathrm{tejF}_{2n}$ of topological even Jacobi forms of index $2n$ [BY26], and also to the truncated spectrum of quasi-modular forms [Dev22], which can be described as $\mathit{tmf} \wedge (S//\nu)^{\leq 4n}$, where $S//\nu$ denotes the $E_1$-cone of the Hopf map $\nu \in \pi_3(S)$, and the superscript denotes the $4n$-skeleton. Similarly, the spectrum $\mathit{tmf} \wedge \mathbb{C}P^n$ is closely related to the spectrum $\mathrm{tjf}_n$ of topological Jacobi forms studied in [BM25]. However, we need detailed knowledge about the life and death of particular unstable classes in $\mathit{tmf}^0(P^n)$, so we cannot directly use these results and instead perform computations with the algebraic Atiyah–Hirzebruch spectral sequence.

In order to lift elements from $\mathit{tmf}^0(P^n)$ to $\pi^0_s(P^n)$, we will use the *tmf*-Hurewicz homomorphism $h : \pi^0_s(P^n) \to \mathit{tmf}^0(P^n)$ mentioned above. Computing the image of this homomorphism for all $n$ would be significantly harder than our computations below, so instead, we use the cofiber sequences

$$P^{n-1} \to P^n \to S^{kn},$$

where $k = 2$ for $\mathbb{C}P^n$ and $k = 4$ for $\mathbb{H}P^n$, in order to use the complete computations of the *tmf*-Hurewicz homomorphisms for spheres

$$h : \pi^0_s(S^{kn}) \to \mathit{tmf}^0(S^{kn})$$

already mentioned above.

## 4. Analysis at $p = 3$

In this section we prove the $p = 3$ part of Theorem A. In Section 4.1 and Section 4.2, we study $\mathit{tmf}^*\mathbb{H}P^n$ and $\mathit{tmf}^*\mathbb{C}P^n$ using the Atiyah–Hirzebruch spectral sequence. In Section 4.3, we use the *tmf*-Hurewicz homomorphism [BS23] together with our Atiyah–Hirzebruch analysis to produce nontrivial elements in $\pi_s^0(\mathbb{C}P^n)$ for certain $n$. In Section 4.4, we deduce geometric consequences.

4.1. **The Atiyah–Hirzebruch spectral sequence for $\mathbb{H}P^n$ at $p = 3$.** At $p = 3$, the ring of derived modular forms is given by

$$\mathrm{Ext}^{**}_{(A,\Gamma)}(A, A) \cong \mathbb{Z}_{(3)}[\alpha, \beta, c_4, c_6, \Delta]/(\alpha^2, c_4^3 + c_6^2 - 1728\Delta)$$

where $|\alpha| = (4, 1)$, $|\beta| = (12, 2)$, $|c_4| = (8, 0)$, $|c_6| = (12, 0)$, and $|\Delta| = (24, 0)$ (see [Bau08, Sec. 5] for details). Since up to suspension, $\mathbb{H}P^2$ is the cone on the Hopf map $\nu \in \pi_3^s$, we have that in $\pi_* tmf[z]$, $d_4(z) = z^2\nu$. Since $\nu$ is represented by $\alpha$ in derived modular forms, we have the fundamental differential $d_4(z) = z^2\alpha$ there.

Figure 1 shows the $E_4$-term of the algebraic Atiyah–Hirzebruch spectral sequence for $\mathbb{H}P^3$, where the differentials $d_4(z) = z^2\alpha$ follows from the above discussion and $d_4(z^2) = 2z^3\alpha$ from the Leibniz rule. Here and in most of the following charts, we omit the ideal generated by $c_4$ and $c_6$, which in any case consists of classes outside the Hurewicz image.

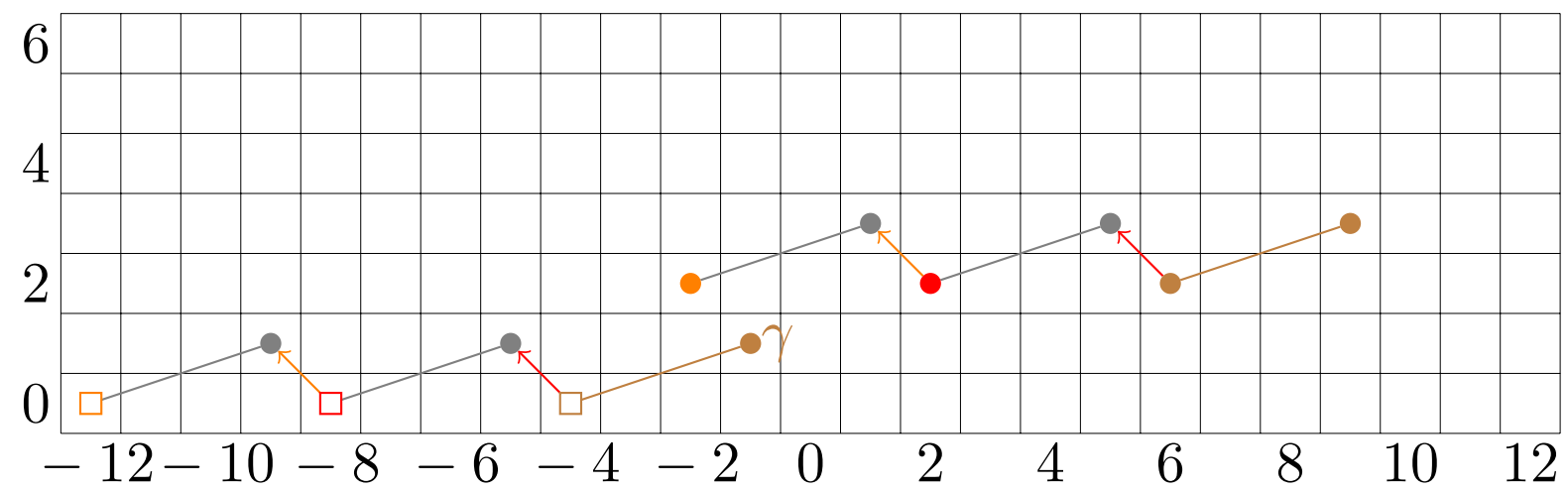

FIGURE 1. The algebraic Atiyah–Hirzebruch spectral sequence for $\mathbb{H}P^3$ at $p = 3$

Now in $\pi_*(S)$, the Toda bracket equation $\beta = \langle \alpha, \alpha, \alpha \rangle$ holds and thus we have that

$$z^3\beta = z^3\langle \alpha, \alpha, \alpha \rangle = \langle z^3, \alpha, \alpha \rangle \alpha$$

must be zero since $\langle z^3, \alpha, \alpha \rangle = 0$. Thus there must be a differential

$$d_8(\gamma) = z^3\beta$$

The spectral sequence collapses at $E_9$ for dimensional reasons since all surviving classes are in Adams–Novikov filtration 0. Also the classical Atiyah–Hirzebruch spectral sequence must therefore collapse at $E_9$.

Since the truncated 3-cell projective spaces $\mathbb{H}P^{3k+3}_{3k+1}$ are homotopy equivalent after localization at $p = 3$, we find that the complete algebraic Atiyah–Hirzebruch spectral sequence for $\mathbb{H}P^\infty$ collapses at $E_9$ as well. Looking more closely at dimension 0, we find the following classes in the $E_4$-term:

$$T_{i,j} = z^{5i+6j}\Delta^j\beta^{2i}$$

When $3 \mid 5i + 6j$, i.e. $3 \mid i$, the class $T_{i,j}$ is a $d^8$-boundary and thus represents 0 in $tmf^*(\mathbb{H}P^{5i+6j})$. In the other cases,

$$d_4(T_{i,j}) = zT_{i,j}\alpha$$

Thus $T_{i,j}$ represents a permanent cycle in the algebraic Atiyah–Hirzebruch spectral sequence for $\mathbb{H}P^{5i+6j}$ (since $Tz = 0$ in this case), and it extends to $\mathbb{H}P^{5i+6j+1}$ if and only if $d_4(T_{i,j}) = 0$ in the *tmf*-based Atiyah–Hirzebruch spectral sequence. Since $\beta^2\alpha = 0$ in $tmf_*$, we see that this is the case whenever $j > 0$, and the classes $T_{i,0}$ are in filtration 0 and hence not in the image of the Hurewicz map.

We conclude that this study at $p = 3$ did not give rise to candidates for free $S^3$-actions, but the computations will be useful to look further for $S^1$-actions.

4.2. **The Atiyah–Hirzebruch spectral sequence for $\mathbb{C}P^n$ at $p = 3$.** Let $\rho$ denote the standard 4-dimensional representation of $S^3$ and consider the cofiber sequence of $S^3$-spaces

$$(S^3/S^1)_+ \wedge S((n+1)\rho) \simeq S^3 \wedge_{S^1} S((n+1)\rho) \to S^{(n+1)\rho} \to S^{\mathrm{ad}} \wedge S((n+1)\rho),$$

where ad, the 3-dimensional adjoint representation of $S^3$, satisfies that $S^{\mathrm{ad}}$ is the unreduced suspension of $S^3/S^1$ as an $S^3$-space. Taking $S^3$-homotopy orbits, we obtain a cofiber sequence

$$\mathbb{C}P^{2n} \to \mathbb{H}P^n \to (\mathbb{H}P^n)^{\mathrm{ad}},$$

where the last space is the Thom space of the adjoint representation, and the connecting map

$$(\mathbb{H}P^n)^{\mathrm{ad}} \to \Sigma\mathbb{C}P^{2n} = (\mathbb{C}P^{2n})^{\mathrm{ad}}$$

is the degree-shifting transfer. Stably, the cofiber sequence

$$(\mathbb{H}P^n)^{\mathrm{ad}-1} \to \mathbb{C}P^{2n} \to \mathbb{H}P^n$$

splits when localized at $p = 3$. Indeed, consider the fibration $f : BS^1 \to BSU(2)$ with fiber $S^2$. The Becker–Gottlieb transfer $t : BSU(2) \to BS^1$ splits this fibration stably after localizing at 3 because in homology, $f \circ t$ is multiplication by the Euler characteristic of the fiber $\chi(S^2) = 2$. The above splitting follows by restricting to $4n$-skeleta. In particular, at $p = 3$,

$$tmf^*(\mathbb{C}P^{2n}) \cong tmf^*(\mathbb{H}P^n) \oplus tmf^*((\mathbb{H}P^n)^{\mathrm{ad}-1}),$$

and the same splitting holds true for cohomotopy groups.

We study the algebraic Atiyah–Hirzebruch spectral sequence for $tmf^*((\mathbb{H}P^n)^{\mathrm{ad}-1})$, as we did before for $\mathbb{H}P^n$. Let $z \in A_{-2}(D\mathbb{C}P^n)$ be the generator of the 2-cell (note that the degree of $z$ in this section is different from the degree of the element called $z$ in the previous section) and denote the corresponding generators of $(\mathbb{H}P^n)^{\mathrm{ad}-1}$ by $z^{2i+1}$. The algebraic Atiyah–Hirzebruch spectral sequence for $(\mathbb{H}P^n)^{\mathrm{ad}-1}$ is a module over the one for $\mathbb{H}P^n$, and we have that $d_4(z) = z^3\alpha$. Thus, the $E_4$-term of the spectral sequence for $(\mathbb{H}P^2)^{\mathrm{ad}-1}$ is as displayed in Figure 2.

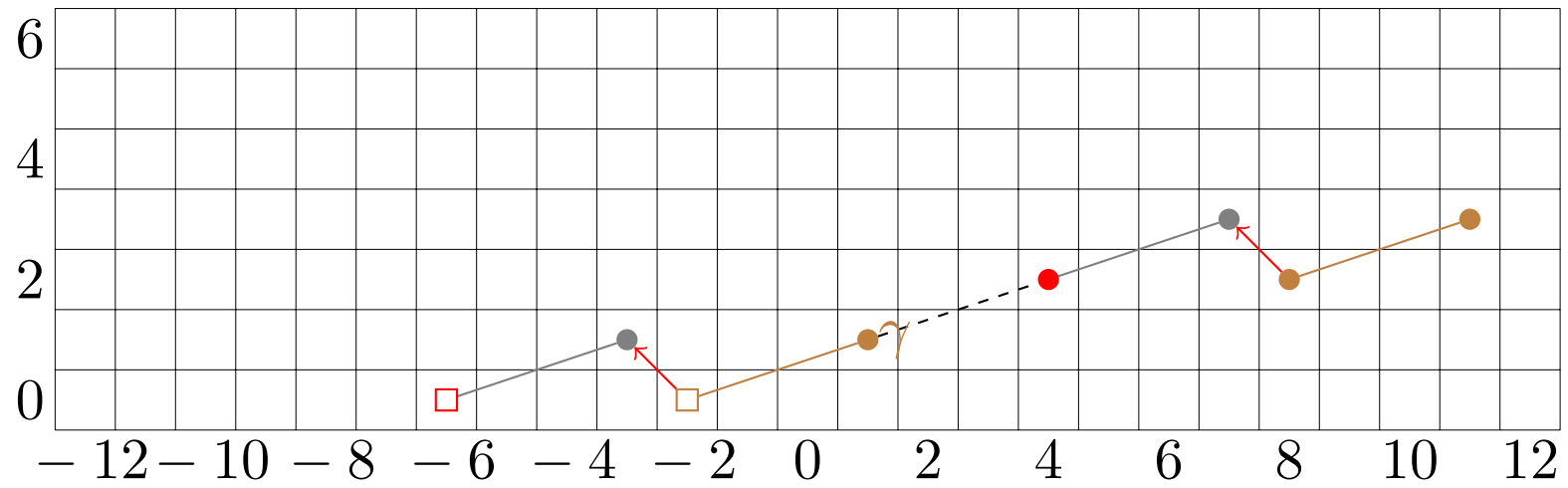

FIGURE 2. The $E_4$-term of the algebraic Atiyah–Hirzebruch spectral sequence for $(\mathbb{H}P^2)^{\mathrm{ad}-1}$

There is an $\alpha$-extension from $\gamma$ to $z^3\beta$ due to the Toda bracket $\langle \alpha, \alpha, \alpha \rangle = \beta$, and the spectral sequence collapses at $E_5$. The remaining $4i+2$-cells of $(\mathbb{H}P^n)^{\mathrm{ad}-1}$ generate modules isomorphic to $\mathrm{Ext}(A_*(\mathbb{H}P^3))$, up to shifts.

In dimension 0, there are no classes in $\mathit{tmf}^0(\mathbb{H}P^2)^{\mathrm{ad}\,-1}$. The other possible candidates in $\mathit{tmf}^0(\mathbb{H}P^n)^{\mathrm{ad}\,-1}$ thus all come from $\mathit{tmf}^{-2}(\mathbb{H}P^n)$:

$$S_{i,j} = z^{5(2i+1)+12j}\Delta^j\beta^{2i+1}.$$

By the same argument as before, only classes with $i = 0$ can possibly contribute. These are only cycles if furthermore $3|j$, so that we conclude:

**Lemma 4.1.** The classes

$$[z^{36j+5}\Delta^{3j}\beta] \in \mathrm{Ext}_{(A,\Gamma)}(A,A)[z]/(z^{36j+7})$$

are infinite cycles, representing nonzero classes in $\mathit{tmf}^0(\mathbb{C}P^{36j+6})$ that do not extend to $\mathit{tmf}^0(\mathbb{C}P^{36j+7})$, and whose image in $\mathit{tmf}^0(S^{72j+13})$ under the attaching map $S^{72j+13} \to \mathbb{C}P^{36j+6}$ of the top cell of $\mathbb{C}P^{36j+7}$ is $\Delta^{3j}\beta\alpha$.

4.3. **Infinite families in the stable cohomotopy of projective spaces, $p = 3$.** Consider the cofiber sequence for $n$ even

$$\mathbb{C}P^{n-2} \xrightarrow{i} \mathbb{C}P^n \xrightarrow{p} \mathbb{C}P^n_{n-1}.$$

Applying the *tmf*-Hurewicz homomorphism yields a map of long exact sequences

$$\begin{array}{ccccccccc} \cdots & \longleftarrow & \pi^0_s\mathbb{C}P^{n-1} & \longleftarrow & \pi^0_s\mathbb{C}P^n & \xleftarrow{p^*} & \pi^0_s\mathbb{C}P^n_{n-1} & \longleftarrow & \cdots \\ & & \big\downarrow & & \big\downarrow & & \big\downarrow{\scriptstyle h} & & \\ \cdots & \longleftarrow & \mathit{tmf}^0\mathbb{C}P^{n-1} & \longleftarrow & \mathit{tmf}^0\mathbb{C}P^n & \xleftarrow{p^*} & \mathit{tmf}^0\mathbb{C}P^n_{n-1} & \longleftarrow & \cdots \end{array} \tag{5}$$

where the right-most vertical map $h$ is the *tmf*-Hurewicz homomorphisms for the $\mathbb{C}P^n_{n-1}$. Since $\mathbb{C}P^n_{n-1} \simeq S^{2n} \vee S^{2n-2}$ after localizing at 3, the right-most vertical map is just a sum of 3-local Hurewicz homomorphisms for these spheres.

Taking $n = 6 + 36j$, Lemma 4.1 shows that

$$p^*(\Delta^{3j}\beta) = [z^{36j+5}\Delta^{3j}\beta] \neq 0 \in \mathit{tmf}^0\mathbb{C}P^{36j+6}.$$

On the other hand, [BS23, Thm. 6.5] implies that $h(\beta_{9t+1}) = \Delta^{6t}\beta$ and $h(\beta_{9t+6/3}) = \Delta^{6t+3}\beta$ for all $t \geq 0$, where $\beta_{9t+1}$ and $\beta_{9t+6/3}$ are certain nonzero classes in the 3-local stable homotopy groups of spheres. Commutativity of (5) then implies that $p^*(\beta_{9t+1}) \neq 0$ and $p^*(\beta_{9t+6/3}) \neq 0$ for all $t \geq 0$.

Reinterpreting these results using the Atiyah–Hirzebruch spectral sequence, we obtain the following.

**Proposition 4.2.** For all $t \geq 0$, the classes

$$z^{72t+5}\beta_{9t+1} \quad \text{and} \quad z^{72t+41}\beta_{9t+6/3}$$

survive in the Atiyah–Hirzebruch spectral sequences converging to $\pi^*_s\mathbb{C}P^{72t+6}$ and $\pi^*_s\mathbb{C}P^{72t+42}$, respectively.

4.4. **Geometric consequence.** We are now prepared to deduce some geometric consequences following the strategy summarized in Section 2.5.

Since $n = 6 + 36j \equiv 2 \mod 4$, the surgery obstruction vanishes on all the classes in $\pi^0_s\mathbb{C}P^n$, so all of the classes obtained in Proposition 4.2 lift to nontrivial classes in $hS(\mathbb{C}P^n)$. It remains to determine their images under $\sigma$, i.e., to determine on which exotic spheres these classes detect free $S^1$-actions.

As above, let $h : \pi_*^s \to tmf_*$ denote the $tmf$-Hurewicz homomorphism for the sphere. We have $h(\alpha_1) = \alpha$. Since $h$ is a ring homomorphism and $\Delta^{3j}\beta\alpha \neq 0 \in tmf_*$ for all $j \geq 0$, this implies that $\alpha_1\beta_{9t+1} \neq 0$ and $\alpha_1\beta_{9t+6/3} \neq 0 \in \pi_*^s$ for all $t \geq 0$.

This relation implies that there are nontrivial $d_4$-differentials

$$d_4(z^{72t+5}\beta_{9t+1}) = z^{72t+7}\alpha_1\beta_{9t+1} \quad \text{and}$$
$$d_4(z^{72t+41}\beta_{9t+6/3}) = z^{72t+43}\alpha_1\beta_{9t+6/3}$$

in the Atiyah–Hirzebruch spectral sequences converging to $\pi_s^*\mathbb{C}P^{7+72t}$ and $\pi_s^*\mathbb{C}P^{43+72t}$.

Putting everything together, we obtain the following.

**Theorem 4.3.** For all $t \geq 0$, there exist exotic spheres with Pontryagin–Thom invariants $\alpha_1\beta_{9t+1}$ and $\alpha_1\beta_{9t+6/3}$ which admit smooth free $S^1$-actions.

**Remark 4.4.** The very first element above, $\alpha_1\beta_1$, recovers a result of Brumfiel [Bru71, Thm. I.10(iv)].

## 5. Analysis at $p = 2$

In this section, we prove the $p = 2$ part of Theorem A. In Section 5.1, we combine Bredon's results with information from the 2-primary $tmf$-Hurewicz homomorphism for spheres to produce some free $S^1$- and $S^3$-actions. Then, in Section 5.2, we analyze $tmf^*\mathbb{H}P^n$ using the Atiyah–Hirzebruch spectral sequence (AHSS), and in particular the algebraic Atiyah–Hirzebruch spectral sequence (AAHSS), and in Section 5.3, we use the $tmf$-Hurewicz homomorphism [BMQ23] to obtain nontrivial elements in $\pi_s^0(\mathbb{H}P^n)$ for certain $n$. These elements are then used in Section 5.4 to study free $S^3$-actions on exotic spheres.

### 5.1. Examples from Bredon's construction at $p = 2$.

Recall from Proposition 2.17 that we may produce a smooth free $S^1$-action (respectively $S^3$-action) on an exotic $(4\ell-1)$-spheres if its 2-primary Kervaire–Milnor invariant is divisible by $\eta$ (respectively $[\ell]\nu$, where $[\ell]$ is the congruence class modulo 8). We will use the $tmf$-Hurewicz homomorphism

$$h : \pi_*^s \to tmf_*.$$

Observe that if $a, \tilde{x}, \tilde{y} \in \pi_*^s$ and $x, y \in tmf_*$ satisfy $ax = y$, $h(\tilde{x}) = x$, and $h(\tilde{y}) = y$, then $a\tilde{x} = \tilde{y} + z$ for some $z \in \pi_*^s$ with $h(z) = 0$. In particular, if we have a relation of the form $\eta x \neq 0$ in $tmf_*$ for some $x \in tmf_{4k-2}$, then we deduce that there is a non-bp $(4k-1)$-sphere which admits a smooth free $S^1$-action, and similarly, a relation of the form $[k]\nu x \neq 0$ in $tmf_*$ for some $x \in tmf_{4k-4}$ implies the existence of a non-bp $(4k-1)$-sphere admitting a smooth free $S^3$-action.

The following is immediate from [BMQ23].

**Lemma 5.1.** The following relations hold in $tmf_*$:

(1) $\eta \cdot \Delta^{4k}\kappa = \Delta^{4k}\kappa\eta \in tmf_{15+96k}$;
(2) $\nu \cdot \Delta^{4k}q = \Delta^{4k}q\nu \in tmf_{35+96k}$.

All of these elements are in the image of the 2-primary Hurewicz homomorphism.

**Corollary 5.2.** For all $k \geq 0$, there exists an exotic sphere with Pontryagin–Thom invariants lifting $\Delta^{4k}\kappa\eta \in tmf_{15+96k}$ which admits a smooth free $S^1$-action, and there exists an exotic sphere with Pontryagin–Thom invariant a lift of $\Delta^{4k}q\nu \in tmf_{35+96k}$ which admits a smooth free $S^3$-action.

5.2. **The Atiyah–Hirzebruch spectral sequence for $\mathbb{H}P^n$ at $p=2$.** The (A)AHSS is somewhat easier than that for complex projective spaces. If $z \in A^*(\mathbb{H}P^8)$ denotes the generator of the bottom cell of $\mathbb{H}P^8$ then by the Leibniz rule and the fact that $\mathbb{H}P^2$ is the cone on $\nu$, we obtain the differentials

$$d^4(z^i) = iz^{i+1}h_2.$$

We will see that the spectral sequence collapses at $E^{16}$ and splits into a sum of two isomorphic spectral sequences for $\mathbb{H}P^4$ and $\mathbb{H}P^8_5$, respectively. Since we need information for all $\mathbb{H}P^n$, we instead study the long exact sequences induced by the functor $\mathrm{Ext}_{(A,\Gamma)}(A_*, A_*(-))$ for the skeletal cofibrations

$$\mathbb{H}P^{i-1} \to \mathbb{H}P^i \to S^{4i},$$

considering it a two-stage spectral sequence.

Figure 3 shows the ring of derived modular forms $\mathrm{Ext}_{(A,\Gamma)}(A,A)$; the reader is advised to refer to this chart for the naming of the indecomposable elements. In this and all other charts, a square denotes an infinite cyclic group, a dot denotes a copy of $\mathbb{Z}/2\mathbb{Z}$, a circle around a class denotes a nontrivial extension by $\mathbb{Z}/2\mathbb{Z}$ (so that a circled dot represents a copy of $\mathbb{Z}/4\mathbb{Z}$), and a diamond denotes a copy of $\mathbb{Z}/2\mathbb{Z}[h_1]$, an infinite $h_1$-tower. In all remaining charts, we will not display the classes in the ideal $(c_4, c_6, \Delta)$ in the zero line to avoid clutter. All of these are infinite cycles in the AAHSS for any $\mathbb{H}P^n$.

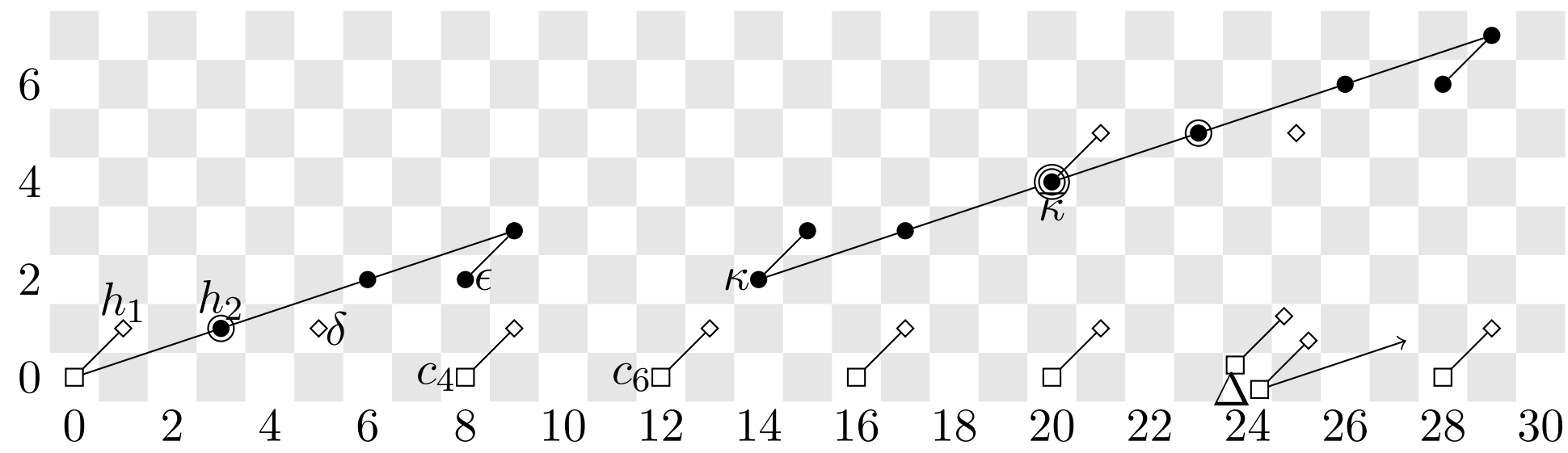

FIGURE 3. The ring $\mathrm{Ext}_{(A,\Gamma)}(A,A)$ of derived modular forms

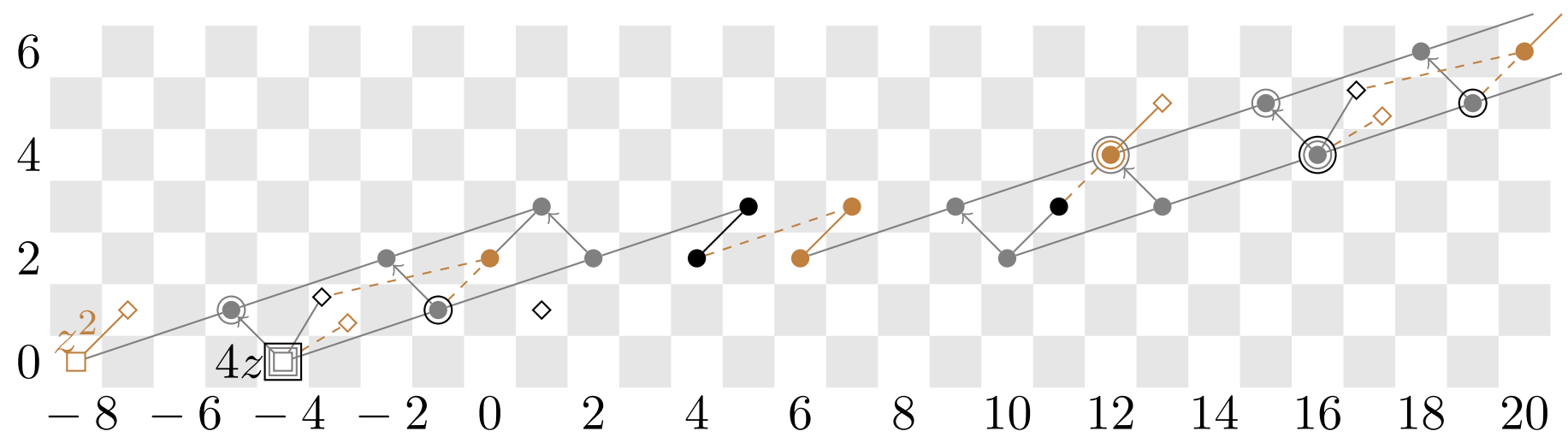

FIGURE 4. The long exact sequence in Ext for $\mathbb{H}P^1 \to \mathbb{H}P^2$

5.2.1. $\mathbb{H}P^2$. The long exact sequence computing $tmf^*(\mathbb{H}P^2)$ is shown in Figure 4. There are several multiplicative extensions present in the target. As is customary, denote by the symbol $\{xy\}$ a class in the $E_\infty$-term which is represented in $E_1$ by $xy$, but where at least one of $x$ and $y$ do not survive to the $E_\infty$-term (and hence $\{xy\}$ is not a product there).

$$
\begin{aligned}
\{4z\}h_1 = \langle z^2, h_2, 4\rangle h_1 = z^2\langle h_2, 4, h_1\rangle \quad &= z^2\delta \\
\{2zh_2\}h_1 = \langle z^2, h_2, 2h_2\rangle h_1 = z^2\langle h_2, 2h_2, h_1\rangle \quad &= z^2\epsilon \\
\{z\epsilon\}h_2 = \langle z^2, h_2, \epsilon\rangle h_2 = z^2\langle h_2, \epsilon, h_2\rangle \quad &= z^2\kappa h_1 \\
\{z\kappa h_1\}h_1 = \langle z^2, h_2, \kappa h_1\rangle h_1 = z^2\langle h_2, \kappa h_1, h_1\rangle \quad &= 2z^2\overline{\kappa} \\
\{zh_1\}h_2 = \langle z^2, h_2, h_1\rangle h_2 = z^2\langle h_2, h_1, h_2\rangle \quad &= z^2\epsilon.
\end{aligned}
$$

The abutment is shown in Figure 5, along with the Adams–Novikov $d^3$-differentials.

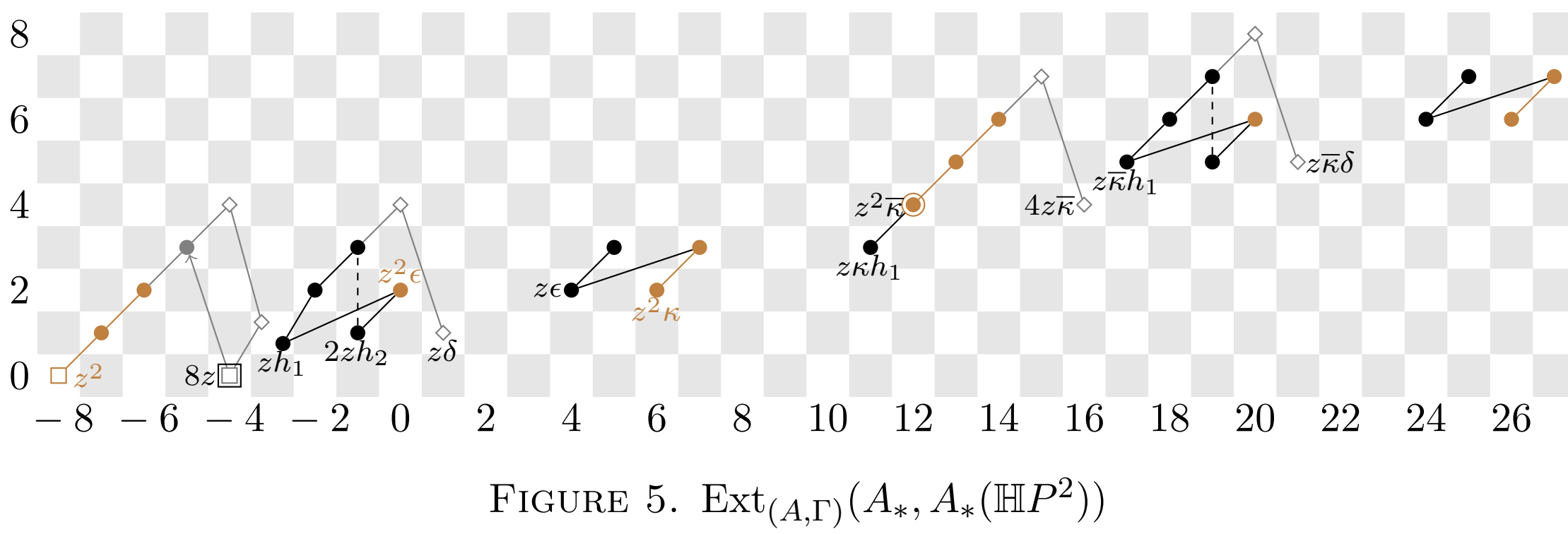

Figure 5. $\mathrm{Ext}_{(A,\Gamma)}(A_*, A_*(\mathbb{H}P^2))$

This computation, along with all higher Adams–Novikov differentials, also appeared in [Tom26]. The reader may also compare with the Adams spectral sequence computations of Bruner and Rognes [BR21, Ch. 8, Sec. 12.3].

5.2.2. $\mathbb{H}P^3$. The long exact sequence for

$$\mathbb{H}P^2 \to \mathbb{H}P^3 \to S^{12}$$

is displayed in Figure 6, and the $E^3$-term of the descent spectral sequence in Figure 7.

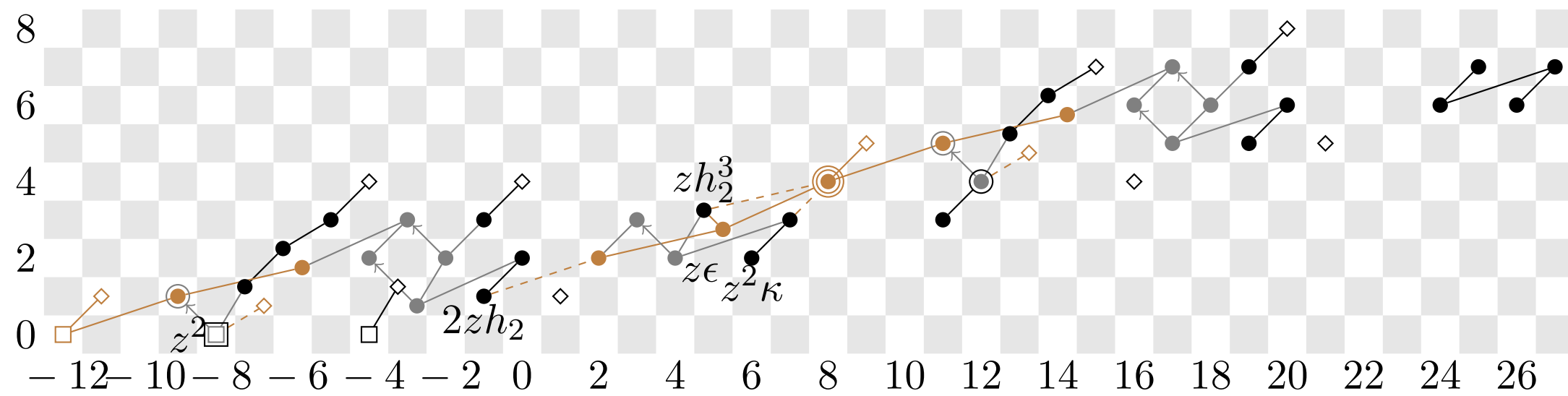

Figure 6. The long exact sequence in $\mathrm{Ext}_{(A,\Gamma)}(A_*, A_*(-)$ coming from the cofiber sequence $\mathbb{H}P^2 \to \mathbb{H}P^3 \to S^{12}$

The Leibniz rule says that $d(z^2) = 2zd(z) = 2z^3h_2$. The Massey product $\{zh_1\} = \langle z^2, h_2, h_1\rangle$ along with the Massey product Leibniz rule implies

$$d(\{zh_1\}) = \langle d(z^2), h_2, h_1\rangle = \langle 2z^3h_2, h_2, h_1\rangle = z^3\langle 2h_2, h_2, h_1\rangle = z^3\epsilon.$$

This actually corresponds to a $d^8$-differential in the AAHSS, since the class $\{zh_1\}$ is defined on the 4-cell and $z^3\epsilon$ on the 12-cell.

Similarly, we have

$$d(\{z\epsilon\}) = d(\langle z^2, h_2, \epsilon\rangle) = \langle 2z^3h_2, h_2, \epsilon\rangle = z^3\langle 2h_2, h_2, \epsilon\rangle = z^3\kappa h_1,$$

another $d^8$-differential in the AAHSS. There is only one more class that for dimensional considerations might support a nontrivial differential, namely, $z^2\kappa$. However, $d(z^2\kappa) = 2z^3\kappa h_2 = 0$.

We study the multplicative extensions by 2, $h_1$, and $h_2$:

$$\begin{aligned}
\{2z^2\}h_1 &= \langle z^3, 2h_2, 2\rangle h_1 = z^3\langle 2h_2, 2, h_1\rangle && = z^3\delta \\
\{2zh_2\}h_2 &= \langle z^3, 2h_2, h_2, 2h_2\rangle h_2 = z^3\langle 2h_2, h_2, 2h_2, h_2\rangle && = z^3\kappa \\
\{z^2\kappa h_1\}h_1 &= \langle z^3, 2h_2, \kappa h_1\rangle h_1 && = 4z^3\bar{\kappa} \\
2\{zh_2^3\} &= \{2zh_2\}h_2^2 && = z^3\kappa h_2 \\
\{zh_2^3\}h_2 &= && 2z^3\bar{\kappa}.
\end{aligned}$$

The last equation holds because twice it holds by the penultimate equation, and so it has to hold up to a summand of $4z^3\bar{\kappa}$, which we may ignore by changing the generator in that degree.

There cannot be an extension of $2\{z^2h_1^2\} = 2\{z^2h_1\}h_1$ because $2h_1 = 0$, and no other extensions by 2, $h_1$, or $h_2$ are possible for dimensional reasons.

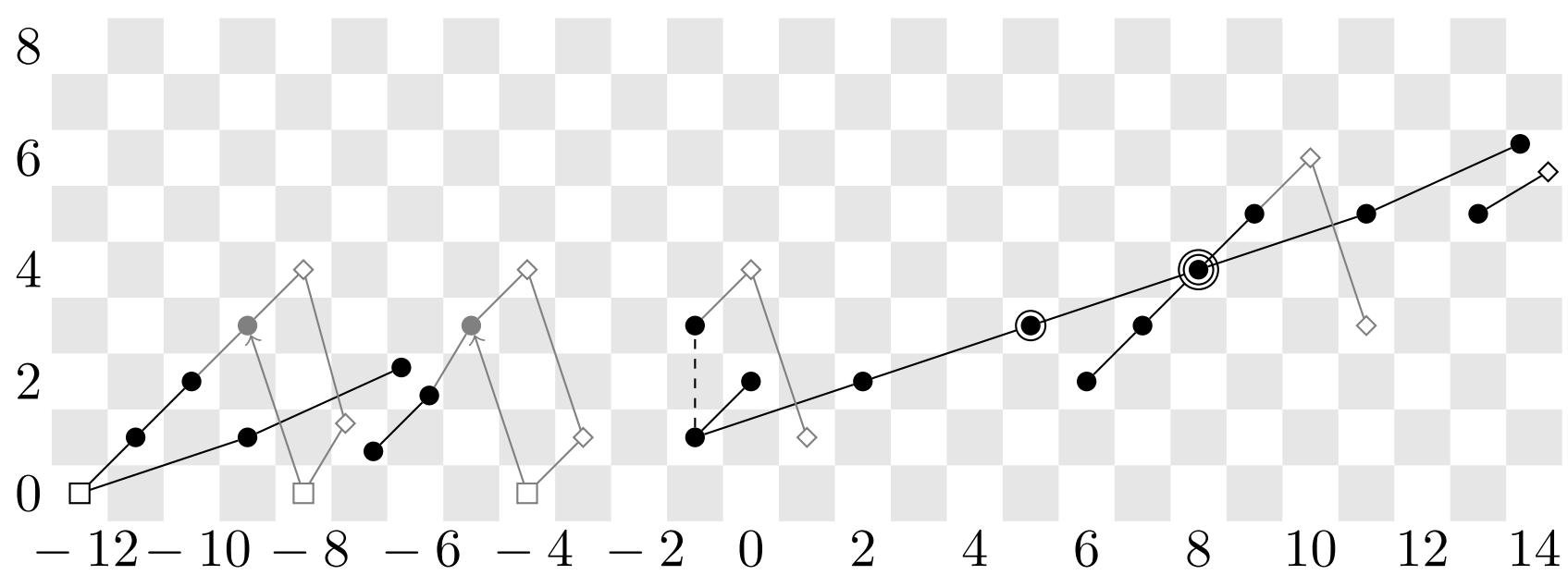

FIGURE 7. $\mathrm{Ext}_{(A,\Gamma)}(A_*, A_*(\mathbb{H}P^3))$ with $d^3$-differentials

5.2.3. $\mathbb{H}P^4$. As usual, we use the previous Ext computation and the cofiber sequence

$$\mathbb{H}P^3 \to \mathbb{H}P^4 \to S^{16}.$$

The corresponding long exact sequence is displayed in Figure 8, and the $E^3$-term of the descent spectral sequence in Figure 9.

By the Leibniz rule, we have $d(z^3) = 3z^4h_2$. We also have

$$d(\{z^2h_1\}) = d(\langle z^3, 2h_2, h_1\rangle) = \langle z^4 3h_2, 2h_2, h_1\rangle = z^4\epsilon,$$

another instance of a $d^8$-differential, and all other differentials are implied by multiplicativity, including the $d^{12}$-differential

$$d(\{2zh_2\}) = z^4\kappa,$$

which follows from $\{2zh_2\}h_2 = z^3\kappa$.

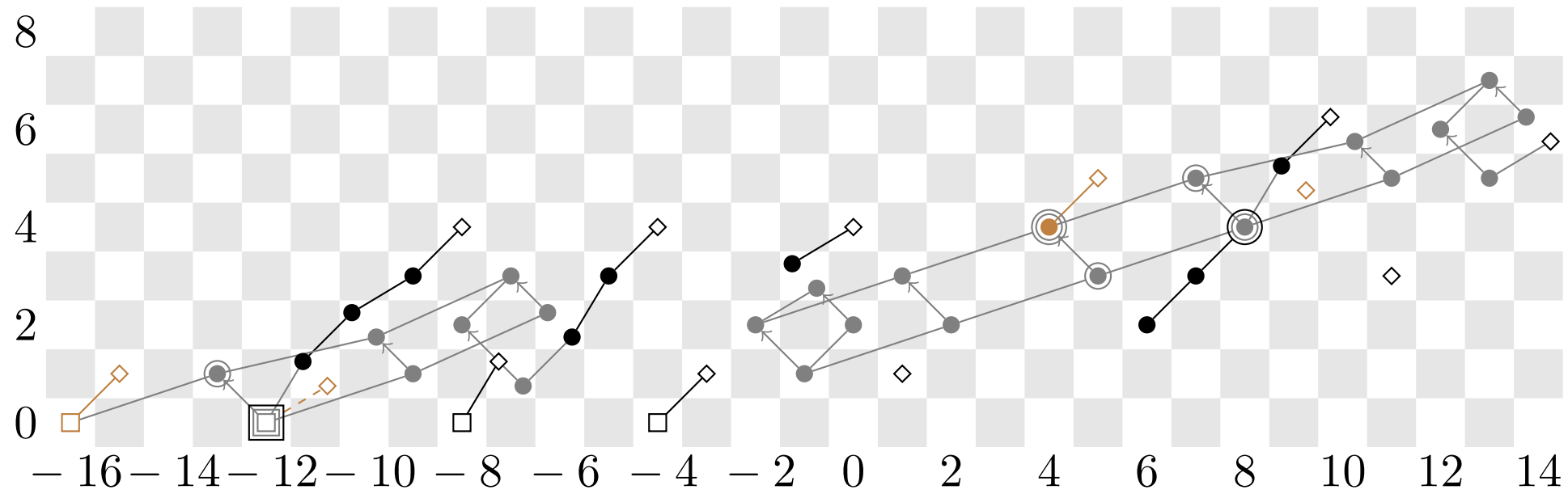

FIGURE 8. The long exact sequence in $\mathrm{Ext}_{(A,\Gamma)}(A_*, A_*(-)$ coming from the cofiber sequence $\mathbb{H}P^3 \to \mathbb{H}P^4 \to S^{16}$

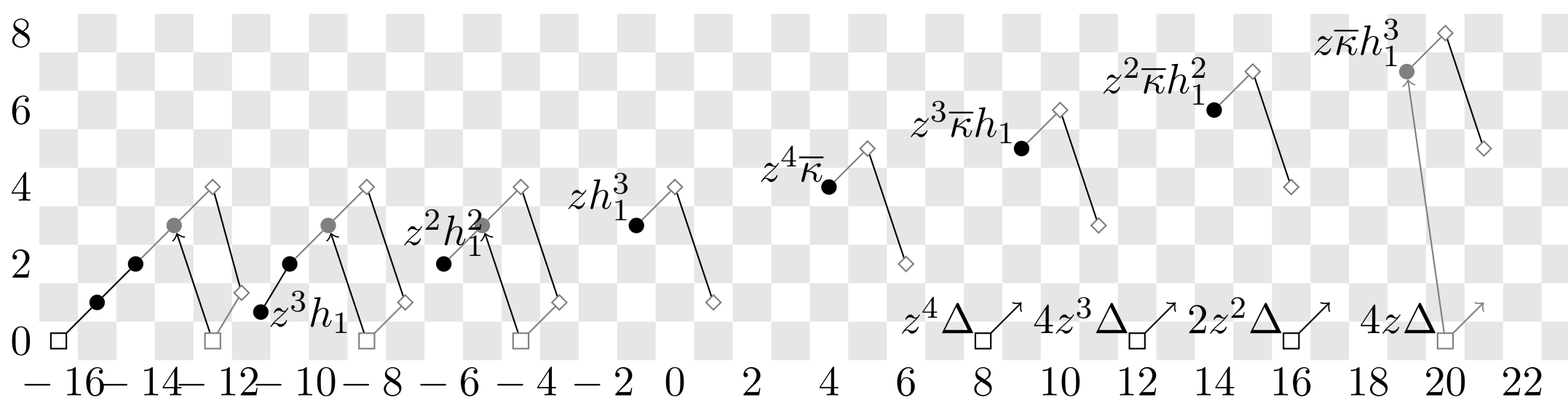

FIGURE 9. $\mathrm{Ext}_{(A,\Gamma)}(A_*, A_*(\mathbb{H}P^4))$

We see that there is no room whatsoever for extensions by 2, $h_1$, or $h_2$.

We observe the following properties:

- The above computations of $\mathrm{Ext}^{s,t}_{(A,\Gamma)}(A, A^*(X))$ for $X = \mathbb{H}P^4$ are also valid for $X = \mathbb{H}P^8_5$ since they only depended on the attaching maps of $4i+4$-cells to $4i$-cells, and $4h_2 = 0$.
- $\mathrm{Ext}^{s,t}_{(A,\Gamma)}(A, A^*(\mathbb{H}P^4)) = 0$ unless $4 \mid (t-s)-s = t-2s$. Since all AAHSS differentials have $(s,t)$-bidegree $(1,0)$, it follows that the AAHSS for $\mathbb{H}P^{4n+k}$ $(0 \leq k \leq 3)$ is equal to a direct sum of $n$ copies of the AAHSS for $\mathbb{H}P^4$, shifted accordingly, and one copy of the AAHSS for $\mathbb{H}P^k$. In particular, the AAHSS collapses at $E^{13}$.
- The descent spectral sequence for $\mathbb{H}P^8$ has a differential $d^3(z^4) = \{z^5h_1^3\}$. At $E_4$, it is concentrated in bidegrees $(s,t)$ with $s \leq 2$ and $t - s \not\equiv 3 \pmod 4$.

The last property implies:

**Proposition 5.3.** In the cofiber sequence $\mathbb{H}P^{8n} \to \mathbb{H}P^{8n+k} \to \mathbb{H}P^{8n+k}_{8n+1}$, the connecting map $\Sigma^{-1}\mathbb{H}P^{8n+k}_{8n+1} \to \mathbb{H}P^{8n}$ induces 0 in $tmf_{(2)}$-homology. In particular, we have a splitting

$$tmf \wedge \mathbb{H}P^{8n+k} \simeq \left( \bigoplus_{i=0}^{n-1} tmf \wedge \mathbb{H}P^{8i+8}_{8i+1} \right) \oplus tmf \wedge \mathbb{H}P^{8n+k}_{8n+1}.$$

*Proof.* Let $X$ be a connective spectrum with cells only in dimensions congruent to $-1$ modulo 4 and $f\colon X \to \mathbb{H}P^{8n}$. Suppose inductively that the restriction $f^{(n)}\colon X^{(4n-1)} \to X \to tmf \wedge (\mathbb{H}P^{8n})$ to the $4n-1$-skeleton of $X$ is nullhomotopic. Then $f^{(n+1)}$ factors through $X^{(4n+3)}/X^{(4n-1)} \cong \bigoplus_i S^{4n+3}$. Since $tmf_{4n+3}(\mathbb{H}P^{8i}) = 0$ for all $n$, $f^{(n+1)}$ is nullhomotopic. □

**Remark 5.4.** The computations also suggest that the truncated spectra $\mathbb{H}P^{8n+8}_{8n+1}$ all become equivalent after smashing with *tmf*. This is true and can be proved using a Thom isomorphism derived from the fact that the spherical fibration associated to the virtual bundle $8L-1$ is *tmf*-orientable, where $L$ is the canonical quaternionic line bundle on $\mathbb{H}P^\infty$. As this statement is not required for proving our results, we omit the details here.

Due to the 16-periodicity of $\mathrm{Ext}_{(A,\Gamma)}(A, A^*(\mathbb{H}P^\infty))$ and the 32-periodicity of $tmf_*(\mathbb{H}P^\infty)$, these computations enable us to find classes $x \in tmf^0(\mathbb{H}P^i)$ that do not extend to $tmf^0(\mathbb{H}P^{i+1})$, and determine their (nontrivial) image $s(x) \in tmf^0(S^{4i+3})$ under the attaching map $S^{4i+3} \to \mathbb{H}P^i$ of the top cell of $\mathbb{H}P^{i+1}$. We are looking for classes $x \in \mathrm{Ext}_{(A,\Gamma)}(A, A^*(\mathbb{H}P^i))$ that are permanent cycles in the descent spectral sequence, and distinguish two types:

(1) We have that $s(x) \neq 0 \in \mathrm{Ext}_{(A,\Gamma)}(A, A^*(S^{4i+3}))$, and $s(x)$ is not an eventual boundary in the descent spectral sequence of the sphere.
(2) We have $s(x) = 0 \in \mathrm{Ext}_{(A,\Gamma)}(A, A^*(S^{4i+3}))$, supporting an Adams–Novikov differential $d^i(x) = y$ in the descent spectral sequence for $\mathbb{H}P^{i+1}$, where $y$ represents a nonzero class $s([x]) \in tmf^*(S^{4i+3})$.

We creatively call these classes candidates of type 1 and 2, respectively.

5.2.4. *Candidates of type 1.* Since the periodicity of *tmf* is 192 and $16 \mid 192$, we are looking for such classes in dimension $16i$. The ring of derived modular forms has the two periodicity elements $\bar\kappa$ in dimension 20 and $\Delta$ in dimension 24, and any class of dimension greater or equal to 20 is divisible by $\bar\kappa$ or $\Delta$. All of our candidate classes therefore are of the form $z^k\Delta^j\bar\kappa^i x$ with $x$ in a dimension divisible by 4.

For any $k$, we consider the classes $nz^k\Delta^j\bar\kappa^i$ for $4 \nmid n$ of dimension $24j + 20i - 4k$. For $\Delta^j\bar\kappa^i$ to represent a nonzero class in $\pi_* tmf$, we must have $j$ even and

$$i \leq \begin{cases} 5; & j \equiv 0 \pmod 8 \\ 1; & j \equiv 2, 6 \pmod 8 \\ 3; & j \equiv 4 \pmod 8 \end{cases} \tag{6}$$

Furthermore, we have that $d_4(nz^k\Delta^j\bar\kappa^i) = nkz^{k+1}\Delta^j\bar\kappa^i h_2$. Since $\bar\kappa h_2$ is a boundary in the Adams–Novikov spectral sequence, we must have $i = 0$. That leaves us with classes of the form $nz^k\Delta^j$ for $4 \nmid nk$, none of which are in the Hurewicz image.

For $k \equiv 1 \pmod 4$, we also have:

- The classes $\Delta^j z^k h_1^2$. However, $d^4(z^k h_1^2) = z^{k_1} h_1^5$, and $h_1^5$ represents $\eta^5 = 0$.

For $k \equiv 2 \pmod 4$, we have:

- The classes $\Delta^j z^{k-1}\bar\kappa^i\epsilon$.

This class hits $z^{k+1}\kappa h_1$ when the dimensional congruence $6j + 5i \equiv -1 \pmod 4$ is satisfied. Since $\bar\kappa^2\epsilon = 0$, we must have $i = 1$, i.e. $\Delta^j z^{6j+7}\bar\kappa\epsilon$, which hits $\Delta^j z^{6j+9}\bar\kappa\kappa h_1$. The classes $\Delta j\bar\kappa\epsilon$ are only permanent cycles for $j \equiv 0, 1 \pmod 4$ in the Adams–Novikov spectral sequence, thus we can further restrict to classes of the form $\Delta^{4j} z^{24j+7}\bar\kappa\epsilon$ or $\Delta^{4j+1} z^{24j+13}\bar\kappa\epsilon$. Only the latter class has the exponent of $z$ is congruent to 1 mod 4, and we have found our first infinite family of classes:

**Lemma 5.5.** For each $j \geq 1$, the class $z^{24j+13}\Delta^{4j+1}\bar{\kappa}\epsilon \in tmf^0(\mathbb{H}P^{24j+14})$ is nonzero, and under the attaching map of the top cell $S^{96j+59}$ of $\mathbb{H}P^{24j+15}$, it maps to the nontrivial class $\Delta^{4j+1}\bar{\kappa}\kappa h_1$.

For $k \equiv 3 \pmod 4$, we have

- The classes $\Delta^j z^{k-1}\bar{\kappa}^i\epsilon$.

If the dimension is to be 0, we must have $i$ even and thus, since $\bar{\kappa}^2\epsilon = 0$, $i = 0$. The classes $\Delta^j z^{6j+2}\epsilon$ are thus candidates. However, as we argued before, $\Delta^j\epsilon$ is not an Adams–Novikov cycle unless $j \equiv 0, 1 \pmod 4$, thus we need to consider $\Delta 4jz^{24j+2}\epsilon$ and $\Delta^{4j+4}z^{24j+8}$. In this case, only the first case satisfies that the power of $z$ is congruent to 2 modulo 4, and we obtain:

**Lemma 5.6.** For each $j \geq 1$, the class $z^{24j+2}\Delta^{4j}\epsilon \in tmf^0(\mathbb{H}P^{24j+3})$ is nonzero, and under the attaching map of the top cell $S^{96j+15}$ of $\mathbb{H}P^{24j+4}$, it maps to the nontrivial class $\Delta^{4j}\kappa h_1$.

For $k \equiv 0 \pmod 4$, obtain no further classes.

5.2.5. *Candidates of type 2.* For any $k$, we have the classes $\{\Delta^j\bar{\kappa}^i nz^k\}$, where $n = \frac{4}{\gcd(4,n)}$, and which support an Adams–Novikov differential $d^3(\{nz^k\}) = z^{k+1}h_1^3$ unless $8 \mid k$. Since $\bar{\kappa}^i\eta^3 = 0$ in $\pi_* tmf$, only the case $i = 0$ is of interest. However, in that case, none of these classes lie in the Hurewicz image. All other $d^3$-differentials hit classes that represent zero in *tmf*-homology. Tominaga's computations [Tom26] of the Adams–Novikov spectral sequence for $\mathbb{H}P^2$ reveal no candidates of type 2 there. It would be interesting to see if there are such candidates for $\mathbb{H}P^3, \ldots, \mathbb{H}P^8$, but we won't do this here.

5.3. **Infinite families in the stable cohomotopy of projective spaces,** $p = 2$**.** For all $n$, we have a cofiber sequence

$$\mathbb{H}P^{n-1} \xrightarrow{i} \mathbb{H}P^n \to S^{4n}.$$

Applying the *tmf*-Hurewicz homomorphism yields a map of long exact sequences

$$\begin{array}{ccccccccc} \cdots & \longleftarrow & \pi_s^0\mathbb{H}P^{n-1} & \longleftarrow & \pi_s^0\mathbb{H}P^n & \xleftarrow{p^*} & \pi_s^0 S^{4n} & \longleftarrow & \cdots \\ & & \downarrow & & \downarrow & & \downarrow h & & \\ \cdots & \longleftarrow & tmf^0\mathbb{H}P^{n-1} & \longleftarrow & tmf^0\mathbb{H}P^n & \xleftarrow{p^*} & tmf^0 S^{4n} & \longleftarrow & \cdots \end{array} \tag{7}$$

where the right-most vertical map $h$ is the *tmf*-Hurewicz homomorphism for $S^{4n}$.

Taking $n = 24j + 13$ and $n = 24j + 2$, we have

$$p^*(\Delta^{4j+1}\bar{\kappa}\epsilon) = [z^{24j+13}\Delta^{4j+1}\bar{\kappa}\epsilon] \neq 0 \in tmf^0(\mathbb{H}P^{24j+13}),$$

$$p^*(\Delta^{4j}\epsilon) = [z^{24j+2}\Delta^{4j}\epsilon] \neq 0 \in tmf^0(\mathbb{H}P^{24j+2}),$$

by Lemma 5.5 and Lemma 5.6, respectively. On the other hand, [BMQ23, Thm. 1.2] implies that $\Delta^{4j+1}\bar{\kappa}\epsilon$ and $\Delta^{4j}\epsilon$ admit lifts along the *tmf*-Hurewicz homomorphism; let us write $v_2^{16j}q\bar{\kappa}$ and $v_2^{16j}\epsilon$, respectively, for these lifts. Commutativity of (7) implies that

$$[z^{24j+13}v_2^{16j}q\bar{\kappa}] := p^*(v_2^{16j}q\bar{\kappa}) \neq 0 \in \pi_s^0\mathbb{H}P^{24j+13},$$

$$[z^{24j+2}v_2^{16j}\epsilon] := p^*(v_2^{16j}\epsilon) \neq 0 \in \pi_s^0\mathbb{H}P^{24j+2}$$

for all $j \geq 0$. Reinterpreting in terms of the AHSS, we have:

**Proposition 5.7.** For all $j \geq 0$, the classes

$$z^{24j+13}v_2^{16j}q\bar{\kappa} \quad \text{and } q^{24j+2}v_2^{16j}\epsilon$$

survive in the Atiyah–Hirzebruch spectral sequences converging to $\pi_s^*\mathbb{H}P^{24j+13}$ and $\pi_s^*\mathbb{H}P^{24j+2}$, respectively.

5.4. **Geometric consequences.** We now follow the strategy of Section 2.5 to deduce the existence of smooth free $S^3$-actions on certain exotic spheres.

The surgery obstruction vanishes on all classes in $\pi_s^0\mathbb{H}P^n$ for all $n$, so all of the classes of Proposition 5.7 lift to nontrivial classes in $hS(\mathbb{H}P^n)$. We wish to determine on which homotopy spheres these classes detect free $S^3$-actions.

There are two possibilities, depending on relations in the stable stems which we do not attempt to determine here. To state the theorem, for all $j \geq 0$, let us write $v_2^{16j+4}\bar{\kappa}\kappa\eta$ and $v_2^{16j}\kappa\eta$ for lifts of $\Delta^{4j+1}\bar{\kappa}\kappa h_1$ and $\Delta^{4j}\kappa\eta$ along the *tmf*-Hurewicz homomorphism for the sphere. Such lifts exist by [BMQ23, Thm. 1.2].

**Theorem 5.8.** Let $j \geq 0$.

(1) If $\nu \cdot v_2^{16j}q\bar{\kappa} \neq 0 \in \pi^s_{55+96j}$, then there exists an exotic sphere with Pontryagin–Thom invariant $\nu \cdot v_2^{16j}q\bar{\kappa}$ which admits a smooth free $S^3$-action.
(2) If $\nu \cdot v_2^{16j}q\bar{\kappa} = 0 \in \pi^s_{55+96j}$, then there exists an exotic sphere with Pontryagin–Thom invariant $v_2^{16j+4}\bar{\kappa}\kappa\eta$ which admits a smooth free $S^3$-action.

*Proof.* Fix $j \geq 0$. Suppose that $\nu \cdot v_2^{16j}q\bar{\kappa} \neq 0$ and consider the AHSS converging to $\pi_s^*\mathbb{H}P^{24j+14}$. By Proposition 5.7, the class $z^{24j+13}v_2^{16j}q\bar{\kappa}$ is not the target of a differential. Since $\nu \cdot v_2^{16j}q\bar{\kappa} \neq 0 \in \pi^s_{55+96j}$, there is a nontrivial $d_4$-differential

$$d_4(z^{24j+13}v_2^{16j}q\bar{\kappa}) = z^{24j+14}\nu \cdot v_2^{16j}q\bar{\kappa} \neq 0.$$

It follows that $\sigma([z^{24j+13}v_2^{16j}q\bar{\kappa}]) = \nu \cdot v_2^{16j}q\bar{\kappa}$. This proves part (1).

Now suppose that $\nu \cdot v_2^{16j}q\bar{\kappa} = 0$ and consider the AHSS converging to $\pi_s^*\mathbb{H}P^{24j+15}$. As above, we know the class $z^{24j+13}v_2^{16j}q\bar{\kappa}$ is not the target of a differential, but now we have

$$d_4(z^{24j+13}v_2^{16j}q\bar{\kappa}) = z^{24j+14}\nu \cdot v_2^{16j}q\bar{\kappa} = 0$$

so $z^{24j+13}v_2^{16j}q\bar{\kappa}$ survives to $E_8$.

We now claim that $z^{24j+15}v_2^{20j}\bar{\kappa}\eta$ also survives to $E_8$. It is a permanent cycle for degree reasons, so it suffices to show it is not the target of a $d_4$-differential. Assume towards a contradiction that it were. Then $v_2^{20j}\bar{\kappa}\eta$ would be divisible by $2\nu$, and thus its *tmf*-Hurewicz image $\Delta^{4j+1}\bar{\kappa}\kappa h_1$ would equal $2\nu$ times an element in the *tmf*-Hurewicz image. But this is not the case, as can be verified from inspection of [Bau08] and [BMQ23, Thm. 1.2]. Thus $z^{24j+15}v_2^{20j}\bar{\kappa}\eta$ survives to $E_8$.

In summary, we have shown that the lifts of the source and target of the $d_8$-differential (Lemma 5.5)

$$d_8(z^{24j+13}\Delta^{4j+1}\bar{\kappa}\epsilon) = z^{24j+15}\Delta^{4j+1}\bar{\kappa}\kappa h_1$$

in the AHSS converging to $tmf^*\mathbb{H}P^{24j+15}$ survive to the $E_8$-page of the AHSS converging to $\pi^*\mathbb{H}P^{24j+15}$. Naturality then implies

$$d_8(z^{24j+13}v_2^{16j}q\epsilon) = z^{24j+15}v_2^{20j}\bar{\kappa}\kappa h_1$$

in the latter spectral sequence, which then implies $\sigma([z^{24j+13}v_2^{16j}q\epsilon]) = v_2^{20j}\bar{\kappa}\kappa h_1$. The completes the proof of part (2). □

A similar argument proves the following.

**Theorem 5.9.** Let $j \geq 0$.

(1) If $2\nu \cdot v_2^{16j}\epsilon \neq 0 \in \pi^s_{11+96j}$, then there exists an exotic sphere with Pontryagin–Thom invariant $2\nu \cdot v_2^{16j}\epsilon$ which admits a smooth free $S^3$-action.
(2) If $2\nu \cdot v_2^{16j}\epsilon = 0 \in \pi^s_{11+96j}$, then there exists an exotic sphere with Pontryagin–Thom invariant $v_2^{16j}\kappa\eta$ which admits a smooth free $S^3$-action.

**Remark 5.10.** When $j = 0$, inspection of [IWX23] shows that $\nu \cdot q\bar{\kappa} = 0$ and $2\nu \cdot \epsilon = 0$, so we obtain smooth free $S^3$-actions on exotic spheres with Pontryagin–Thom invariants $v_2^4\bar{\kappa}\kappa\eta$ and $\kappa\eta$. The latter example was also proven recently by Basu and Kasilingam [BK22] using explicit low-dimensional computations.

KTH Royal Institute of Technology
*Email address*: `tilmanb@kth.se`

University of Virginia
*Email address*: `mbp6pj@virginia.edu`